\documentclass[amstex,12pt, amssymb]{article}

\usepackage{mathtext}
\usepackage[cp1251]{inputenc}
\usepackage[T2A]{fontenc}
\usepackage[dvips]{graphicx}
\usepackage{amsmath}
\usepackage{amssymb}
\usepackage{amsxtra}
\usepackage{latexsym}
\usepackage{ifthen}

\newcounter{lemma}[section]

\newcounter{corollary}[section]

\newcounter{remark}[section]

\newcounter{theorem}[section]

\newcounter{proposition}[section]

\numberwithin{equation}{section}

\def\Xint#1{\mathchoice
   {\XXint\displaystyle\textstyle{#1}}%
   {\XXint\textstyle\scriptstyle{#1}}%
   {\XXint\scriptstyle\scriptscriptstyle{#1}}%
   {\XXint\scriptscriptstyle\scriptscriptstyle{#1}}%
   \!\int}
\def\XXint#1#2#3{{\setbox0=\hbox{$#1{#2#3}{\int}$}
     \vcenter{\hbox{$#2#3$}}\kern-.5\wd0}}
\def\dashint{\Xint-}

\def\cc{\setcounter{equation}{0}
\setcounter{figure}{0}\setcounter{table}{0}}

\begin{document}

\markboth{\centerline{V. RYAZANOV, R. SALIMOV AND E. SEVOST'YANOV}}
{\centerline{CONVERGENCE AND COMPACTNESS OF THE SOBOLEV MAPPINGS}}

\author{{V. RYAZANOV, R. SALIMOV AND E. SEVOST'YANOV}\\}

\title{
{\bf CONVERGENCE AND COMPACTNESS \\ OF THE SOBOLEV MAPPINGS}}

\maketitle

\large \begin{abstract} First of all, we establish compactness of
continuous mappings of the Orlicz--Sobolev classes
$W^{1,\varphi}_{\rm loc}$ with the Calderon type condition on
$\varphi$ and, in particular, of the Sobolev classes $W^{1,p}_{\rm
loc}$ for $p>n-1$ in ${\Bbb R}^n\,,$ $n\ge 3\,,$ with one fixed
point. Then we give a series of theorems on convergence of the
Orlicz--Sobolev homeomorphisms and on semicontinuity in the mean of
dilatations of the Sobolev homeomorphisms. These results lead us to
closeness of the corresponding classes of homeomorpisms. On this
basis, we come to criteria of compactness of classes of Sobolev's
homeomorphisms with two fixed points. Finally, we show the precision
of the found conditions.
\end{abstract}


\large \cc
\section{Introduction}\label{1}

The present paper is a natural continuation of our preceding works
\cite{KRSS} and \cite{RS$_3$}. In the last work \cite{RS$_3$}, we
gave foundations of the convergence theory for general
ho\-me\-o\-mor\-phisms in space and, on this basis, developed the
compactness theory for the so--called ring
$Q$--ho\-me\-o\-mor\-phisms that take an important part in the
mapping theory as well as in the theory of the Beltrami equations,
see e.g. the monographs \cite{GRSY} and \cite{MRSY}, and the papers
\cite{KPR} and \cite{KRSS}. Here we give some results of our study
of the corresponding convergence and compactness problems for the
Orlicz--Sobolev mappings and, in particular, for the Sobolev
homeo\-morphisms based on the mentioned theory of ring
$Q$-homeomorphisms that was motivated by the ring definition of
Gehring for quasiconformal mappings, see e.g. \cite{Ge}.


First of all, recall the minimal definitions related to the Sobolev
spaces $W^{1,p}$, $p\in [1, \infty ) $. Given an open set $U$ in
$\mathbb{R}^n$, $n\ge 2,$ $C_0^{ \infty }(U)$ denotes the collection
of all functions $\varphi : U \to \mathbb{R}$ with compact support
having continuous partial derivatives of any order. Now, let $u$ and
$v: U \to {\Bbb R}$ be locally integrable functions. The function
$v$ is called the {\bf distributional
derivative}\index{distributional derivative} $u_{x_i}$ of $u$ in the
variable $x_i$, $i=1,2,\ldots , n$, $x=(x_1,x_2,\ldots , x_n)$, if
\begin{equation}\label{eqSTR2.21}
\int\limits_{U} u\, \varphi_{x_i} \,dm(x)=- \int\limits_{U} v\,
\varphi \ dm(x) \ \quad \forall \ \varphi \in C_{0}^{\infty}(U)\ .
\end{equation} \medskip
Here $dm(x)$ corresponds to the Lebesgue measure in $\mathbb{R}^n$.
The {\bf Sobolev classes} $W^{1,p}(U)$ consist of all functions $u:
U \to {\Bbb R}$ in $L^p(U)$ with all distributional derivatives of
the first order in $L^p(U)$. A function $u: U \to {\Bbb R}$ belongs
to $W^{1,p}_{\mathrm{loc}}(U)$ if $u \in W^{1,p}(U_*)$ for every
open set $U_*$ with a compact closure in $ U.$  We use the
abreviation $W^{1,p}_{\mathrm{loc}}$ if $U$ is either defined by the
context or not essential. The similar notion is introduced for
vector-functions $f: U \to \mathbb{R}^m$ in the component-wise
sense.
\medskip

The concept of the distributional derivative was introduced by
Sobolev, see \cite{So}. It is known that a continuous function $f$
belongs to $W^{1,p}_{\rm loc}$ if and only if $f\in ACL^{p}$, i.e.,
if $f$ is locally absolutely continuous on a.e. straight line which
is parallel to a coordinate axis and if all the first partial
derivatives of $f$ are locally integrable with the power $p$, see
e.g. 1.1.3 in \cite{Maz}.

\medskip

In what follows, $D$ is a domain in a finite-dimensional Euclidean
space. Following Orlicz, see \cite{Or1}, given a convex increasing
function $\varphi:[0,\infty)$ $\rightarrow[0,\infty)$,
$\varphi(0)=0$, denote by $L^{\varphi}$ the space of all functions
$f:D\rightarrow{\Bbb R}$ such that
\begin{equation}\label{eqOS1.1}\int\limits_{D}\varphi\left(\frac{|f(x)|}
{\lambda}\right)\,dm(x)<\infty\end{equation} for some $\lambda>0$.
$L_{\varphi}$ is called the {\bf Orlicz space}.  If
$\varphi(t)=t^p$, then we write also $L^{p}$. In other words,
$L^{\varphi}$ is the cone over the class of all functions
$g:D\rightarrow{\Bbb R}$ such that
\begin{equation}\label{eqOS1.2}
\int\limits_{D}\varphi\left(|g(x)|\right)\,dm(x)<\infty\end{equation}
which is also called the {\bf Orlicz class}, see \cite{BO}.

\medskip

The {\bf Orlicz-Sobolev class} $W^{1,\varphi}_{\rm loc}(D)$ is the
class of locally integrable functions $f$ given in $D$ with the
first distributional derivatives whose gradient $\nabla f$ has a
modulus $|\nabla f|$ that belongs locally in $D$ to the Orlicz
class. Note that by definition $W^{1,\varphi}_{\rm loc}\subseteq
W^{1,1}_{\rm loc}$. Later on, we also write $f\in W^{1,\varphi}_{\rm
loc}$ for a locally integrable vector-function $f=(f_1,\ldots,f_m)$
of $n$ real variables $x_1,\ldots,x_n$ if $f_i\in W^{1,1}_{\rm loc}$
and
\begin{equation}\label{eqOS1.2a}
\int\limits_{D}\varphi\left(|\nabla
f(x)|\right)\,dm(x)<\infty\end{equation} where $|\nabla
f(x)|=\sqrt{\sum\limits_{i,j}\left(\frac{\partial f_i}{\partial
x_j}\right)^2}$. Note that in this paper we use the notation
$W^{1,\varphi}_{\rm loc}$ for more general functions $\varphi$ than
in the classical Orlicz classes sometimes giving up the condition on
convexity of $\varphi$. Note also that the Orlicz--Sobolev classes
are intensively studied in various aspects at present, see e.g.
\cite{KRSS} and further references therein.

\medskip

Now, let $f$ be a continuous map of a domain $D$ in ${\Bbb R}^n$,
$n\geqslant2$, into ${\Bbb R}^n$. If $f$ has all the first partial
derivatives at a point $x\in D$, then we denote by $\Vert
f'(x)\Vert$ the matrix norm of the Jacobian matrix $f'$ of $f$ at
$x$, i.e., $||f'(x)||=\sup\limits_{h\in{\Bbb R}^n,|h|=1}|f'(x)\cdot
h|$, and by $J_f(x)$ its Jacobian, i.e., ${\rm det}f'(x)$. The {\bf
outer dilatation} of $f$ at such a point $x$ is the quantity
\begin{equation}\label{eqOS1.3}
K_f(x)=\frac{\Vert f'(x)\Vert^n}{|J_f(x)|}
\end{equation}
if $J_f(x)\ne0$, $K_f(x)=1$ if $f'(x)=0$, and $K_f(x)=\infty$ at the
rest points in\-clu\-ding points without first partial derivatives.
Further we also use the dilatation
\begin{equation}\label{eqOS1.3.1}
P_{f}\left( x\right) \,=\,{K^{\frac{1}{n-1}}_{f}(x)}\,.
\end{equation}

Finally, recall that a homeomorphism $f$ between domains $D$ and
$D^{\prime}$ in ${\Bbb R}^n$, $n\geqslant2$, is called of {\bf
finite distortion} if $f\in W^{1,1}_{\rm loc}$, $J_f(x)\ge 0$ and
$K_f(x)$ is finite a.e. First this notion was introduced on the
plane for $f\in W^{1,2}_{\rm loc}$ in the work \cite{IS}. Later on,
this condition was changed by $f\in W^{1,1}_{\rm loc}$ in ${\Bbb
R}^n$, $n\geqslant2$, but with the additional condition $J_f\in
L^1_{\rm loc}$ in the monograph \cite{IM}.

\medskip

Note that the above additional condition $J_f\in L^1_{\rm loc}$ in
the definition of the mappings with finite distortion in \cite{IM}
can be omitted for homeomorphisms. Indeed, for each homeomorphism
$f$ between domains $D$ and $D'$ in ${\Bbb R}^n$ having all the
first partial derivatives a.e. in $D$, there is a set $E$ of the
Lebesgue measure zero such that $f$ satisfies $(N)$-property by
Lusin on $D\setminus E$ and
\begin{equation}\label{eqOS1.1.1}
\int\limits_{A}|J_f(x)|\,dm(x)=|f(A)|\end{equation} for every Borel
set $A\subset D\setminus E$, see e.g. 3.1.4, 3.1.8 and 3.2.5 in
\cite{Fe}.

\medskip

On the basis of (\ref{eqOS1.1.1}), it is easy to prove the following
useful statement.

\medskip

\begin{proposition}\label{prOS2.4} {\sl
Let $f$ be an ACL homeomorphism of a domain $D$ in ${\Bbb R}^n$,
$n\geqslant2$, into ${\Bbb R}^n$. Then
$$
(i) \quad  f\in W^{1,1}_{\mathrm loc} \quad  \mbox{if} \quad P_f\in
L^1_{\mathrm loc}\ ,
$$
$$
(ii) \quad  f\in W^{1,\frac{n}{2}}_{\mathrm loc} \quad  \mbox{if}
\quad K_f\in L^1_{\mathrm loc}\ ,
$$
$$
(iii) \quad  f\in W^{1, n-1}_{\mathrm loc} \quad  \mbox{if} \quad
K_f\in L^{n-1}_{\mathrm loc}\ ,
$$
$$
(iv) \quad  f\in W^{1, p}_{\mathrm loc},\ p>n-1 \quad  \mbox{if}
\quad K_f\in L^{{\gamma}}_{\mathrm loc},\ {\gamma}>n-1\ ,
$$
$$
(v) \quad  f\in W^{1, p}_{\mathrm loc},\ p=n{\gamma}/(1+{\gamma})\ge
1 \quad \mbox{if} \quad K_f\in L^{{\gamma}}_{\mathrm loc},\
{\gamma}\ge 1/(n-1)\ .
$$
These conclusions and the estimates (\ref{eqOSadm.1}) are also valid
for all ACL mappings $f: D\to{\Bbb R}^n$ with $J_f\in L^1_{\mathrm
loc}$.}\end{proposition}

\medskip

\begin{proof} Indeed, by the H\"{o}lder inequality applied on a compact
set $C$ in $D$, we obtain on the basis of (\ref{eqOS1.1.1}) the
following estimates of the first partial derivatives
\begin{equation}\label{eqOSadm.1}
\Vert\partial_i f\Vert_p\leqslant\Vert
f^{\prime}\Vert_p\leqslant\Vert K_f^{1/n}\Vert_s\cdot\Vert
J_f^{1/n}\Vert_n\leqslant\Vert
K_f\Vert^{1/n}_{\gamma}\cdot|f(C)|^{1/n}<\infty \end{equation} if
$K_f\in L^{{\gamma}}_{\mathrm loc}$ for some $\gamma\in(0,\infty )$
because $\Vert f^{\prime}(x)\Vert = K_f^{1/n}(x)\cdot J_f^{1/n}(x)$
a.e. where $\frac{1}{p}=\frac{1}{s}+\frac{1}{n}$ and $s={\gamma}n$,
i.e., $\frac{1}{p}=\frac{1}{n}\left(\frac{1}{{\gamma}}+1\right)$.
\end{proof} $\Box$

\bigskip

We sometimes use the estimate (\ref{eqOSadm.1}) with no comments to
obtain corollaries.

\bigskip

Later on, we also often use the notations $\mathrm I$, $\bar{\mathrm
I}$, $\mathbb R$, $\overline{\mathbb R}$, ${\Bbb R}^+$,
$\overline{{\Bbb R}^+}$ and $\overline{{\Bbb R}^n}$ for
$[1,\infty)$, $[1,\infty]$, $(-\infty,\infty)$, $[-\infty ,\infty
]$, $[0,\infty)$, $[0,\infty ]$ and ${\Bbb R}^n\cup\{\infty\}$,
correspondingly.

\bigskip

Moreover, we denote by $\mathrm B(x,r)$, $x\in{\Bbb R}^n$, $n\ge2$,
$r>0$, the open ball centered at $x$ with the radius $r$, i.e.,
$\mathrm B(x,r)=\{ z\in{\Bbb R}^n:\ |z-x|<r\ \}$, and set ${\Bbb
B}^n=\mathrm B(0,1)$.


\cc
\section{ On Compactness of Orlicz--Sobolev Mappings  }\label{11}

Recall definitions related to normal and compact families of
mappings in metric spaces. Let $(X, d)$ and $\left(X^{\,\prime},
d^{\,\prime}\right)$ be metric spaces with distances $d$ and
$d^{\,\prime},$ respectively. A family $\frak{F}$ of continuous
mappings $f:X\rightarrow X^{\,\prime}$ is said to be {\bf normal} if
every sequence of mappings $f_j\in\frak{F}$ has a subsequence
$f_{j_m}$ converging uniformly on each compact set $C\subset X$ to a
continuous mapping $f$. If in addition $\frak{F}$ is {\bf closed}
with respect to the locally uniform convergence, i.e.,
$f\in\frak{F}$, then the family is called {\bf compact}.

\medskip

Normality is closely related to the following notion. A family
$\frak{F}$ of mappings $f:X\rightarrow X^{\,\prime}$ is said to be
{\bf equicontinuous at a point} $x_0\in X$ if for every
$\varepsilon>0$ there is $\delta>0$ such that
$d^{\,\prime}(f(x),f(x_0))<\varepsilon$ for all $f\in\frak{F}$ and
$x\in X$ with $d(x,x_0)<\delta.$ The family $\frak{F}$ is called
{\bf equicontinuous} if $\frak{F}$ is equicontinuous at every point
$x_0\in X.$

\medskip

Now, let us formulate the fundamental Calderon result in \cite{Ca},
p. 208.

\medskip

\begin{proposition}\label{prOS2.4} {\sl Let $\varphi:{\Bbb R}^+\rightarrow{\Bbb R}^+$
be an increasing functi\-on with $\varphi(0)=0$ and the condition
\begin{equation}\label{eqOS2.5}A\ \colon=\  \int\limits_{0}^{\infty}\left[\frac{t}{\varphi(t)}\right]^
{\frac{1}{k-1}}\ dt\ <\ \infty\end{equation} for a natural number
$k\ge2$ and let $f:D\rightarrow{\Bbb R}$ be a continuous function
given in a domain $D\subset{\Bbb R}^k$ of the class
$W^{1,\varphi}(D)$. Then
\begin{equation}\label{eqOS2.6}{\rm diam}\,
f(C)\le\alpha_{k}\,A^{\frac{k-1}{k}}\left[\
\int\limits_{C}\varphi\left(|\nabla
f|\right)\,dm(x)\right]^{\frac{1}{k}}\end{equation} for every cube
$C\subset D$ whose adges are oriented along coordinate axes where
$\alpha_k$ is a constant depending only on $k$.}\end{proposition}

\bigskip

\begin{remark}\label{remOS2.1} Here it is not essential that the
fuction $\varphi$ is (strictly !) increasing. Indeed, let $\varphi$
is only nondecreasing. Going over, in case of need, to the new
function
$$
\tilde{\varphi}_{\varepsilon}(t):=\varphi(t)+\sum\limits_{i}{\varphi}_i^{(\varepsilon)}(t)
$$
where
$$
{\varphi}_i^{(\varepsilon)}(t):=\varepsilon\,\frac{2^{-i}}{(b_i-a_i)}\int\limits_{0}^{t}\chi_i(t)\,dt
$$
and $\chi_i$ is a numbering of the characteristic functions of the
intervals of constancy $(a_i,b_i)$ of the function $\varphi$, we see
that
$\varphi(t)\leqslant\tilde{\varphi}_{\varepsilon}(t)\leqslant\varphi(t)+{\varepsilon}$
and, thus, the condition (\ref{eqOS1.2a}) on $C$ and the condition
(\ref{eqOS2.5}) hold for the (strictly !) increasing function
$\tilde{\varphi}_{\varepsilon}$. Letting $\varepsilon\to0$, we
obtain the estimate  (\ref{eqOS2.6}) with the initial function
$\varphi$, see e.g. Theorem I.12.1 in \cite{Sa}.

The function $\left(t/{\varphi}(t)\right)^{1/(k-1)}$ can have a
nonintegrable singularity at zero. However, it is clear that the
behavior of the function ${\varphi}$ about zero is not essential for
the estimate (\ref{eqOS2.6}). Indeed, we may apply the estimate
(\ref{eqOS2.6}) with the replacements $A\mapsto A_*$ and
$\varphi\mapsto \varphi_*$ where \begin{equation}\label{eqOS2.55}
A_*:=\left[{\frac{1}{\varphi(t_*)}}\right]^{\frac{1}{k-1}}\ +\
\int\limits_{t_*}^{\infty}\left[\frac{t}{\varphi(t)}\right]^{\frac{1}{k-1}}\,dt<\infty\end{equation}
and $\varphi_*(0)=0$, $\varphi_*(t)\equiv\varphi(t_*)$ for
$t\in(0,t_*)$ and $\varphi_*(t)=\varphi(t)$ for $t\geqslant t_*$ if
$\varphi(t_*)>0$. Hence, in particular, the normalization
$\varphi(0) = 0$ in Proposition \ref{prOS2.4} evidently has no
valuation, too.
\end{remark}

\medskip
Given a domain $D$ in ${\Bbb R}^n,$ $n\ge 2,$ a nondecreasing
function $\varphi:\overline{{\Bbb R}^+}\rightarrow \overline{{\Bbb
R}^+},$ $M\in[0,\infty)$ and $x_0\in D$, denote by
$\frak{F}^{\varphi}_M$ the family of all continuous mappings
$f:D\rightarrow{\Bbb R}^m,$ $m\ge 1,$ of the class $W^{1,1}_{\rm
loc}$ such that $f(x_0)=0$ and
\begin{equation}\label{eqOS8.111}
\int\limits_{D}\varphi\left(|\nabla f|\right)\,dm(x)\le M\ .
\end{equation} We also use the notation
$\frak{F}^{p}_M$ for the case of the function $\varphi(t)=t^p,$
$p\in [1, \infty).$ Set
\begin{equation}\label{eqOS8.111.v}
t_0\ =\ \sup\limits_{\varphi(t)=0}\ t\ ,\ \ \ \ \ \ t_0\ =\ 0\ \ \
\mbox{if}\ \ \ \varphi(t)>0\ \ \  \forall t\in\overline{{\Bbb R}^+}
\end{equation} and
\begin{equation}\label{eqOS8.111.w}
T_0\ =\ \inf\limits_{\varphi(t)=\infty}\ t\ ,\ \ \ \ \ \ T_0\ =\
\infty\ \ \ \mbox{if}\ \ \ \varphi(t)<\infty\ \ \  \forall
t\in\overline{{\Bbb R}^+}\ . \end{equation}

\medskip

In this section we prove the following result, cf. e.g. Theorem 8.1
in \cite{IKO}.

\medskip

\begin{theorem}\label{corCALDERON}
{\sl Let $\varphi:\overline{{\Bbb R}^+}\rightarrow \overline{{\Bbb
R}^+}$ be a nonconstant continuous nondecreasing convex function
such that, for some $\ t_*\in(t_0,\infty)$ and $\
\alpha\in(0,1/(n-1))$,
\begin{equation}\label{eq9.111}
\int\limits_{t_*}^{\infty}\left(\frac{t}{\varphi(t)}\right)^
{\alpha} dt\ <\ \infty\ .\end{equation}   Then
$\frak{F}^{\varphi}_M$ is compact with respect to the locally
uniform convergence in ${\Bbb R}^n$. }
\end{theorem}

\medskip

Here the continuity of the function $\varphi$ is understood in the
sense of the to\-po\-lo\-gy of the extended positive real axis
$\overline{{\Bbb R}^+}$.

\medskip

Recall before its proof that a nondecreasing convex function
$\varphi:\overline{{\Bbb R}^+}\rightarrow \overline{{\Bbb R}^+}$ is
called {\bf strictly convex}, see e.g. \cite{Rud}, if
\begin{equation}\label{eq_1.16}
    \lim\limits_{t \to \infty}\ \frac{\varphi (t)} {t}\ =\ \infty\ .
\end{equation}


\begin{remark}\label{remOS2.1.a} Note that a nonconstant continuous nondecreasing convex function
$\varphi:\overline{{\Bbb R}^+}\rightarrow \overline{{\Bbb R}^+}$
satisfying the condition (\ref{eq9.111}) for some $\alpha>0$ is
strictly convex. Indeed, the slope $\varphi(t)/t$ is a nondecreasing
function if $\varphi$ is convex, see e.g. Proposition I.4.5 in
\cite{Bou}. Hence the condition (\ref{eq9.111}) for $\alpha>0$
implies (\ref{eq_1.16}).
\end{remark}

\medskip

The proof of Theorem \ref{corCALDERON} will be based on the
following lemma.


\begin{lemma}\label{lem7.2.a}
{\sl Let $\varphi:\overline{{\Bbb R}^+}\rightarrow \overline{{\Bbb
R}^+}$ be a nonconstant continuous nondecreasing convex function
with the condition (\ref{eq9.111}) for some $\alpha>0$ and let
$\tilde{\alpha}\in(\alpha,\infty)$. Then $\varphi$ admits the
decomposition $\varphi=\psi\circ\tilde{\varphi}$ where $\psi$ and
$\tilde{\varphi}:\overline{{\Bbb R}^+}\rightarrow \overline{{\Bbb
R}^+}$ are strictly convex and, moreover,
$\tilde{\varphi}\le\varphi$ and $\tilde{\varphi}$ satisfies
(\ref{eq9.111}) with the new $\tilde{\alpha}$.}
\end{lemma}

\medskip

\begin{proof} Note that the convex function $\varphi$ is locally
lipshitz on the interval $(0,T_0)$, where $T_0$ is defined by
(\ref{eqOS8.111.w}), $T_0>t_0$ by continuity and variability of the
function $\varphi$. Consequently, $\varphi$ is locally absolutely
con\-ti\-nuo\-us and, furthermore, differentiable except a countable
collection of points in the given nondegenerate interval and
$\varphi^{\prime}$ is non\-de\-crea\-sing, see e.g. Corollaries 1-2
and Proposition 8 of Section I.4 in \cite{Bou}. Thus, denoting by
$\varphi^{\prime}_+(t)$ the function which coincides with
$\varphi^{\prime}(t)$ at the points of differentiability of
$\varphi$ and $\varphi^{\prime}_+(t)=\lim\limits_{\tau\to
t+0}\varphi^{\prime}(\tau)$ at the rest points in the interval
$[0,T_0)$ and, finally, setting $\varphi^{\prime}_+(t)=\infty$ for
all $t\in [T_0,\infty]$, we have that
\begin{equation}\label{eq_1.16.z}
\varphi(t)\ =\ \varphi(0)\ +\
\int\limits^t_0\varphi^{\prime}_+(\tau) \ d\tau\ \ \ \ \ \  \forall\
t\in\overline{{\Bbb R}^+}\ .\end{equation}

By monotonicity of the function $\varphi^{\prime}_+$, calculating
its averages over the segments $[0,t]$ and $[t/2,t]$,
correspondingly, we obtain from (\ref{eq_1.16.z}) the two-sided
estimate
\begin{equation}\label{eq_1.16.v}
\frac{1}{2}\ \varphi^{\prime}_+\left({t}/{2}\right)\ \le\
\frac{\varphi(t)\ -\ \varphi(0)}{t}\ \le\ \varphi^{\prime}_+(t) \ \
\ \ \ \  \forall\ t\in\overline{{\Bbb R}^+}\ .\end{equation} The
inequalities (\ref{eq_1.16.v}) show that the condition
(\ref{eq9.111}) is equivalent to the following
\begin{equation}\label{eq9.111.z}
I\ \colon =\
\int\limits_{t_*}^{\infty}\frac{dt}{[\varphi^{\prime}_+(t)]^
{\alpha}} \ <\ \infty\ .\end{equation}

Again by monotonicity of $\varphi^{\prime}_+$, the condition
(\ref{eq9.111.z}) implies that $\varphi^{\prime}_+(t)\to\infty$ as
$t\to\infty$. Thus, $T_*=\sup\limits_{\varphi^{\prime}_+(t)<1}t$ is
finite, $T_*\in [t_0,T_0)$. Set $\lambda = {\alpha}/{\alpha_*}\in
(0,1)$.

Consider the functions $\tilde{\varphi}(t)= \int\limits_0^t
h(\tau)\, d\tau$ and $\psi(s)=\varphi(0) + \int\limits_0^s H(r)\,
dr$ where $h(t)=\varphi^{\prime}_+(t)$ for $t\in [0,T_*)$ and
$h(t)=[\varphi^{\prime}_+(t)]^{\lambda}$ for $t\in [T_*,\infty ]$
and $H(s)=1$ for $s\in [0,S_*)$, $S_*=\varphi_*(T_*)$,
$H(s)=[\varphi^{\prime}_+(\tilde{\varphi}^{-1}(s))]^{1-\lambda}$ for
$s\in [S_*,S_0)$, $S_0=\varphi_*(T_0)$, and $H(s)=\infty$ for $s\in
[S_0,\infty ]$.


By the construction, $\tilde{\varphi}(t)\le\varphi(t)$ for all
$t\in\overline{{\Bbb R}^+}$, the functions $\psi$ and
$\tilde{\varphi}$ as well as $\psi\circ\tilde{\varphi}$ are
nondecreasing and convex, see e.g. Proposition 8 of Section I.4 in
\cite{Bou}, and
\begin{equation}\label{eq9.111.v}
\int\limits_{t_*}^{\infty}\frac{dt}{[\tilde{\varphi}^{\prime}_+(t)]^
{\tilde{\alpha}}}\ =\ I\ <\ \infty\ \end{equation} and, thus,
$\tilde{\varphi}$ satisfies (\ref{eq9.111}) with the new
$\tilde{\alpha}$. Moreover, similarly to (\ref{eq_1.16.v})
\begin{equation}\label{eq_1.16.w}
\frac{\psi(s)\ -\ \psi(0)}{s}\ \ge\ \frac{1}{2}\
H\left({s}/{2}\right)\ \  \ \ \ \ \ \ \forall\ s\in\overline{{\Bbb
R}^+}\end{equation} where the right hand side converges to $\infty$
as $s\to\infty$. Thus, $\psi$ is strictly convex.


Finally, simple calculations by the chain rule show that
$$
(\psi\circ\tilde{\varphi})^{\prime}_+(t)\ =\
\psi^{\prime}_+(\tilde{\varphi}(t))\cdot\tilde{\varphi}_+^{\prime}(t)\
=\ \varphi^{\prime}_+(t)\
$$
except a countable collection of points in $\overline{{\Bbb R}^+}$,
$\psi\circ\tilde{\varphi}(0)=\varphi(0)$ and, consequently,
$\psi\circ\tilde{\varphi}\equiv\varphi$ in view of
(\ref{eq_1.16.z}).
\end{proof}$\Box$

\bigskip

And now, let us give the proof of the main result of this section,
Theorem \ref{corCALDERON}.

\bigskip

\begin{proof}
First, let us show that mappings in $\frak{F}^{\varphi}_M$ are
equicontinuous. Indeed, by Lemma \ref{lem7.2.a} $\varphi$ admits the
decomposition $\varphi=\psi\circ\tilde{\varphi}$ where $\psi$ and
$\tilde{\varphi}:\overline{{\Bbb R}^+}\rightarrow \overline{{\Bbb
R}^+}$ are strictly convex and
\begin{equation}\label{eq9.111.z.1}
\int\limits_{t_*}^{\infty}\left(\frac{t}{\tilde{\varphi}(t)}\right)^
{\frac{1}{n-1}} dt\ <\ \infty\end{equation} for some $t_*>t_0$ and,
moreover, $\tilde{\varphi}\le\varphi$ and hence
\begin{equation}\label{eqOS8.111.z}
\int\limits_{D}\tilde{\varphi}\left(|\nabla f|\right)\,dm(x)\le M\ .
\end{equation}
Given $z_0\in D$ and $\delta>0,$ denote by $C(z_0, \delta)$ the
$n$--dimensional open cube centered at the point $z_0$ with edges
which are parallel to coordinate axes  and whose length is equal to
$\delta.$ Fix $\varepsilon>0.$ Since the function $\psi$ is strictly
convex, the integral of $\tilde{\varphi}(|\nabla f|)$ over $C(z_0,
\delta)\subset D$ is arbitrary small at sufficiently small
$\delta>0$ for all $f\in\frak{F}^{\varphi}_M$, see e.g. Theorem
III.3.1.2  in \cite{Rud}. Thus, by Proposition \ref{prOS2.4} and
Remark \ref{remOS2.1} applied to $\tilde{\varphi}$ we have that
 $|f(z)-f(z_0)|<\varepsilon$ for all $z\in
C(z_0, \delta)$ under some $\delta=\delta(\varepsilon)>0.$

\medskip

Now, let us show that a family $\frak{F}^{\varphi}_M$ is uniformly
bounded on compactums. Indeed, let $K$ be a compactum in $D.$ With
no loss of generality we may consider that $K$ is a connected set
containing the point $x_0$ from the definition of
$\frak{F}^{\varphi}_M$, see e.g. Lemma 1 in \cite{Smol}. Let us
cover $K$ by the collection of cubes $C (z, \delta_z),$ $z\in K,$
where $\delta_z$ corresponds to $\varepsilon:=1$ from the first part
of the proof. Since $K$ is compact, we can find a finite number of
cubes $C_i=C(z_i, \delta_{z_i}),$ $i=1,2,\ldots,N$ that cover $K$.
Note that $D_*:=\bigcup\limits_{i=1}^N C_i$ is a subdomain of $D$
because $K$ is a connected set. Consequently, each point $z_*\in K$
can be joined with $x_0$ in $D_*$ by a polygonal curve with ends of
its segments at points $x_0,x_1,\ldots, x_k, z_*$ in the given order
lying in the cubes with numbers $i_1,\ldots, i_k,$ $x_0\in
C(z_{i_1}, \delta_{z_{i_1}}),$ $z_*\in C(z_{i_k}, \delta_{z_{i_k}})$
and $x_l\in C_{i_l}\cap C_{i_{l+1}},$ $l=1,\ldots, k-1,$ $k\le N-1.$
By the triangle inequality we have that
$$|f(z_*)|\le \sum\limits_{l=0}^{k-1}
|f(x_l)-f(x_{l+1})|+|f(x_k)-f(z_*)|\le N\,.$$
Since $N$ depends on a compactum $K$ only, it follows that
$\frak{F}^{\varphi}_M$ is uniformly bounded on compactums and,
consequently, is normal by the Arzela–-Ascoli theorem, see e.g.
IV.6.7 in \cite{DS}.

\medskip

Finally, show that the class $\frak{F}^{\varphi}_M$ is closed. By
Remark \ref{remOS2.1.a} $\varphi$ is strictly convex and by Theorem
III.3.1.2 in \cite{Rud}, for every $\varepsilon>0$, there is
$\delta=\delta(\varepsilon)>0$ such that $\int\limits_{E} |\nabla
f|\,dm(x)\le \varepsilon$ for all $f\in\frak{F}^{\varphi}_M$
whenever $m(E)<\delta.$ Let $f_j\in \frak{F}^{\varphi}_M$ and
$f_j\rightarrow f$ locally uniformly as $j\rightarrow \infty.$ Then
by Lemma 2.1 in \cite{RSY} we have the inclusion $f\in W_{\rm
loc}^{1,1}.$ By Theorem 3.3 in Ch. III, $\S\,3.4$, of the monograph
\cite{Re},
\begin{equation}\label{eq5.2}
\int\limits_{D}\varphi(|\nabla f|)\ dm(x)\ \le\ M\,,
\end{equation}
i.e., $\frak{F}^{\varphi}_M$  is closed. Thus, the class
$\frak{F}^{\varphi}_M$ is compact.
\end{proof}$\Box$

\medskip
\begin{corollary}\label{cor5.1}
{\sl\, The class $\frak{F}^{p}_M$ is compact with respect to the
locally uniform convergence for each $p\in (n, \infty).$}
\end{corollary}

\medskip
\begin{proof} It is easy to verify that the function $\varphi(t)=t^p$
satisfies the hypotheses of Theorem \ref{corCALDERON} for an
arbitrary number $\alpha\in(1/(p-1), 1/(n-1))$. \end{proof}$\Box$

\bigskip

Recall that the problem of equicontinuity of mappings in the classes
$W^{1,p}$ for $p>n$ was investigated in the well--known paper
\cite{BI}, cf. also \cite{IKO$_1$}. However, the condition $p>n$ is
too restrictive for mappings with finite distortion as it was
cleared already in the plane case, see e.g. \cite{BGR}, \cite{Da},
\cite{GRSY} and \cite{MRSY}, although this condition was natural for
quasiconformal mappings, see e.g. \cite{Boj} and \cite{Ge$_{12}$}.
Hence we will go back to this question once more in Section 7.

\cc
\section{Convergence of Orlicz--Sobolev Homeomorphisms}
\label{4}

In \cite{KRSS}, see Corollary 9.3, it was established by us that
every homeomorphism between domains in ${\Bbb R}^n\,,$ $n\ge 3\,,$
of the Orlicz--Sobolev classes $W^{1,\varphi}_{\mathrm loc}$ with
the Calderon type condition on a nonconstant continuous
nondecreasing function $\varphi: \overline{{\Bbb
R}^+}\to\overline{{\Bbb R}^+}$:
\begin{equation}\label{eqOS3.0a}\int\limits_{t_*}^{\infty}\left[\frac{t}{\varphi(t)}\right]^
{\frac{1}{n-2}}dt<\infty\end{equation} for some $t_*\in(t_0,\infty)$
where $t_0=\sup\limits_{\varphi(0)=0}t$ , see also Remark
\ref{remOS2.1}, and, in particular, of the Sobolev classes
$W^{1,p}_{\mathrm loc}$ with $p>n-1$  is a ring $Q_*$-homeomorphism
at every point $x_0\in{D}$ with $Q_*(x)=\left[K_f(x)\right]^{n-1}$.
Thus, combining this fact with the convergence theory for ring
$Q-$homeomorphisms in the last work \cite{RS$_3$}, Section 4, we
come to the corresponding results on convergence of homeomorphisms
in the classes of Sobolev as well as Orlicz--Sobolev further.

\medskip
\begin{lemma}\label{lem7.2}
{\sl\,Let $D$ be a domain in ${\Bbb R}^n,$ $n\ge 3,$ and let $f_j,$
$j=1,2,\ldots,$ be a sequence of ho\-me\-o\-mor\-phisms of $D$ into
${\Bbb R}^n$ in $W^{1,\varphi}_{\mathrm loc}$ with (\ref{eqOS3.0a})
converging locally uniformly to a mapping $f$ with respect to the
spherical metric. Suppose
\begin{equation}\label{eq6.3.24.A}
\int\limits_{\varepsilon<|x-x_0|<\varepsilon_0}
K_f^{n-1}(x)\cdot\psi^n(|x-x_0|) \ dm(x)=o(I^{n}(\varepsilon,
\varepsilon_0))\quad\quad\forall\ x_0\in D
\end{equation}
where $o(I^{n}(\varepsilon, \varepsilon_0))/I^{n}(\varepsilon,
\varepsilon_0)\to 0$ as $\varepsilon\rightarrow 0$ uniformly with
respect to the parameter $j=1,2,\ldots$ for $\varepsilon_0<{\rm
dist}\,(x_0,
\partial D)$ and a measurable function $\psi(t):(0, \varepsilon_0)\rightarrow [0, \infty]$
such that
\begin{equation}\label{eq6.3.25A}
0<I(\varepsilon,
\varepsilon_0):=\int\limits_{\varepsilon}^{\varepsilon_0}\psi(t)\,\,dt
< \infty\quad\quad\forall\,\varepsilon\in(0, \varepsilon_0)\ .
\end{equation}
Then the mapping $f$ is either a constant in $\overline{{\Bbb R}^n}$
or a homeomorphism into ${\Bbb R}^n.$}
\end{lemma}

\begin{remark}\label{remOS2.1.z} In particular, the conclusion of Lemma \ref{lem7.2} holds for
ho\-meo\-mor\-phisms $f_j$ in the classes $W^{1,p}_{\mathrm loc}$
with $p>n-1$ and for $K_{f_j}^{n-1}(x)\le Q(x)$ a.e. in $D$ with a
measurable function $Q: D\rightarrow (0, \infty)$ such that
\begin{equation}\label{eq6.3.24}
\int\limits_{\varepsilon<|x-x_0|<\varepsilon_0}
Q(x)\cdot\psi^n(|x-x_0|) \ dm(x)=o(I^{n}(\varepsilon,
\varepsilon_0))\quad\quad\forall\ x_0\in D\ .
\end{equation}
\end{remark}

\medskip
\begin{theorem}{}\label{th6.6.1A}{\sl\,
Let $D$ be a domain in ${\Bbb R}^n,$ $n\ge 3,$ $f_j,$
$j=1,2,\ldots,$ be a sequence of ho\-me\-o\-mor\-phisms of $D$ into
${\Bbb R}^n$  in $W^{1,\varphi}_{\mathrm loc}$ with (\ref{eqOS3.0a})
and $K_{f_j}^{n-1}(x)\le Q(x)$ a.e. in $D$ where $Q\in$ FMO. If
$f_j\to f$ locally uniformly, then the mapping $f$ is either a
constant in $\overline{{\Bbb R}^n}$ or a homeomorphism into ${\Bbb
R}^n.$}
\end{theorem}

\medskip
\begin{corollary}\label{cor6.6.2A}{\sl\,
In particular, the limit mapping $f$ is either a constant in
$\overline{{\Bbb R}^n}$ or a homeomorphism of $D$ into
${\Bbb R}^n$ 
%
whenever
$$\overline{\lim\limits_{\varepsilon\rightarrow 0}}\ \
 \dashint_{\mathrm B( x_0 ,\varepsilon)} Q(x)\ \ dm(x) <
 \infty\qquad\qquad\forall\,\,\,x_0
\in D 
$$
or  whenever every $x_0 \in D$ is a Lebesgue point of $Q\in
L^{1}_{\mathrm loc}$.}
\end{corollary}

\medskip
\begin{theorem}{}\label{th6.6.5A}{\sl\, Let $D$ be a domain in ${\Bbb
R}^n,$ $n\ge 3,$ $Q:D \rightarrow \bar{\mathrm I}$ a locally
integrable function such that, for some $\varepsilon(x_0)< {\rm
dist}\, (x_0,\partial D)$,
 \begin{equation}\label{eq6.6.6}
\int\limits_{0}^{\varepsilon(x_0)}\frac{dr}{rq_{x_0}^{\frac{1}{n-1}}(r)}=\infty\quad
\quad \forall\,x_0\in D\end{equation}  where $q_{x_0}(r)$ denotes
the average of $Q(x)$ over the sphere $|x-x_0|=r.$ Suppose $f_j,$
$j=1,2,\ldots,$ is a sequence of ho\-me\-o\-mor\-phisms of $D$ into
${\Bbb R}^n$  in $W^{1,\varphi}_{\mathrm loc}$ with (\ref{eqOS3.0a})
and $K_{f_j}^{n-1}(x)\le Q(x)$ a.e. in $D$. If $f_j\to f$ locally
uniformly, then the mapping $f$ is either a constant in
$\overline{{\Bbb R}^n}$ or a homeomorphism into ${\Bbb R}^n.$ }
\end{theorem}

\medskip
\begin{corollary}\label{cor4.1}
{\sl \, In particular, the conclusion of Theorem \ref{th6.6.5A}
holds if $$q_{x_0}(r)=O\left(\log^{n-1}\frac{1}{r}\right)\qquad
\forall \,x_0\in D\,.$$}
\end{corollary}

\medskip
\begin{corollary}\label{cor6.6.7A}{\sl\, Under hypotheses of
Theorem \ref{th6.6.5A}, the limit mapping $f$ is either a constant
in $\overline{{\Bbb R}^n}$ or a homeomorphism into ${\Bbb R}^n$
provided that the function $Q$ has singularities only of the
logarithmic type of the order which is not more than $n-1$ at every
point $x_0\in D.$}
\end{corollary}

\medskip
\begin{theorem}\label{th4.1}{\sl\, Let $D$ be a domain in ${\Bbb
R}^n,$ $n\ge 3,$ $Q:D \rightarrow \bar{\mathrm I}$ a locally
integrable function such that
\begin{equation}\label{eq4.2}
\int\limits_{\varepsilon<|x-x_0|<\varepsilon_0}\frac{Q(x)}{|x-x_0|^n}\,dm(x)=
o\left(\log^n\frac{1}{\varepsilon}\right)\qquad\forall\, x_0\in D
\end{equation}
as $\varepsilon\rightarrow 0$ for some positive number
$\varepsilon_0=\varepsilon(x_0)<{\rm dist}\,(x_0,
\partial D).$ Suppose $f_j,$
$j=1,2,\ldots,$ is a sequence of ho\-me\-o\-mor\-phisms of $D$ into
${\Bbb R}^n$  in $W^{1,\varphi}_{\mathrm loc}$ with (\ref{eqOS3.0a})
and $K_{f_j}^{n-1}(x)\le Q(x)$ a.e. in $D$. If $f_j\to f$ locally
uniformly, then the mapping $f$ is either a constant in
$\overline{{\Bbb R}^n}$ or a homeomorphism into ${\Bbb R}^n.$
}
\end{theorem}

\bigskip

For every nondecreasing function $\Phi:\bar {\mathrm
I}\rightarrow\overline{{\mathbb R}^+}$, the {\bf inverse function}
$\Phi^{-1}:\overline{{\mathbb R}^+}\rightarrow\bar{\mathrm I}$ can
be well defined by setting
$$\Phi^{-1}(\tau)\ =\ \inf\limits_{\Phi(t)\ge \tau}\ t\ .$$
As usual, here $\inf$ is equal to $\infty$ if the set of $t$ in
$\bar{\mathrm I}$  such that $\Phi(t)\ge \tau$ is empty. Note that
the function $\Phi^{-1}$ is nondecreasing, too. Note also that if
$h:\bar{\mathrm I}\rightarrow\bar{\mathrm I}$ is a sense--preserving
homeomorphism and $\varphi :\bar{\mathrm
I}\rightarrow\overline{{\mathbb R}^+}$ is a nondecreasing function,
then
\begin{equation}\label{eq1!!!}
(\varphi\circ h)^{-1}\ =\ h^{-1}\circ\varphi^{-1} \ . \end{equation}

\medskip

\begin{theorem}{}\label{th4.2}{\sl\,Let $D$ be a domain in ${\Bbb
R}^n,$ $n\ge 3,$ $Q:D \rightarrow\bar{\mathrm I}$ a measurable
function and $\Phi:\bar{\mathrm I}\rightarrow\overline{{\mathbb
R}^+}$ a nondecreasing convex function. Suppose that
\begin{equation}\label{eq2!!}
\int\limits_D\Phi\left(Q(x)\right)\frac{dm(x)}{\left(1+|x|^2\right)^n}\
\le\ M<\infty
\end{equation}
and, for some $\delta>\Phi(1),$
\begin{equation}\label{eq3!}
\int\limits_{\delta}^{\infty}
\frac{d\tau}{\tau\left[\Phi^{-1}(\tau)\right]^{\frac{1}{n-1}}}\ =\
\infty
\end{equation}
 Suppose $f_j,$ $j=1,2,\ldots,$ is a
sequence of ho\-me\-o\-mor\-phisms of $D$ into ${\Bbb R}^n$  in
$W^{1,\varphi}_{\mathrm loc}$ with (\ref{eqOS3.0a}) and
$K_{f_j}^{n-1}(x)\le Q(x)$ a.e. in $D$. If $f_j\to f$ locally
uniformly, then $f$ is either a constant in $\overline{{\Bbb R}^n}$
or a homeomorphism into ${\Bbb R}^n.$}
\end{theorem}

\medskip

\begin{remark}\rm\label{rmkCOM3.1.1}
We may assume in Theorem \ref{th4.2} that the function $\Phi(t)$ is
not convex on the whole segment $\bar{\mathrm I}$ but only on the
segment $[t_*,\infty]$ where $t_*=\Phi^{-1}(\delta)$. Indeed, every
non-decreasing function $\Phi :\bar{\mathrm I}\to\overline{{\mathbb
R}^+}$ which is convex on the segment $[t_*,\infty]$ can be replaced
by a non-decreasing convex function $\Phi_*:\bar{\mathrm
I}\to\overline{{\mathbb R}^+}$ minorizing $\Phi$  in the following
way. Set $\Phi_*(t)\equiv 0$ for $t\in [1,{t_*} ]$,
$\Phi(t)=\varphi(t)$ for $t\in[t_*,T_*]$ and $\Phi_*\equiv \Phi(t)$
for $t\in[T_*,\infty]$, where $\tau =\varphi(t)$ is the line passing
through the point $(0,{t_*})$ and touching the graph of the function
$\tau =\Phi(t)$ at a point $(T_*,\Phi(T_*))$, $T_*\in({t_*},\infty
)$.  By the construction we have that $\Phi_*(t)\le \Phi(t)$ for all
$t\in\bar{\mathrm I}$ and $\Phi_*(t)=\Phi(t)$ for all $t\ge T_*$
and, consequently, the conditions (\ref{eq2!!}) and (\ref{eq3!})
hold for $\Phi_*$ under the same $M$ and every $\delta
> 0$.

Furthermore, by the same reasons it is sufficient to assume in
Theorem \ref{th4.2} that the function $\Phi$ is only minorized by a
nondecreasing convex function $\Psi$ on a segment $[T,\infty ]$ such
that
\begin{equation}\label{eq!3!}
\int\limits_{\delta}^{\infty}
\frac{d\tau}{\tau\left[\Psi^{-1}(\tau)\right]^{\frac{1}{n-1}}}\ =\
\infty
\end{equation}
for some $T\in\mathrm I$ and $\delta>\Psi(T).$ Note that the
condition (\ref{eq!3!}) can be written in terms of the function
$\psi(t)=\log \Psi(t)\ $:
\begin{equation}\label{eq3!!!33}
\int\limits_{\Delta}^{\infty}\ \psi(t)\ \frac{dt}{t^{n^{\prime}}}\
=\ \infty
\end{equation}
for some $\Delta>t_0\in [T,\infty ]$ where
$t_0:=\sup\limits_{\psi(t)=-\infty} t,$ $t_0=T$ if $\psi(T)>-\infty
,$ and where $\frac{1}{n'}+\frac{1}{n}=1$, i.e., $n'=2$ for $n=2$,
$n'$ is decreasing in $n$ and $n'=n/(n-1)\to1$ as $n\to\infty$, see
Proposition 2.3 in \cite{RS$_2$}. It is clear that if the function
$\psi$ is nondecreasing and convex, then the function
$\Phi=e^{\psi}$ is so but the inverse conlusion generally speaking
is not true. However, the conclusion of Theorem \ref{th4.2} is valid
if $\psi^m(t)$, $t\in [T,\infty ]$, is convex and (\ref{eq3!!!33})
holds for $\psi^m$ under some $m\in{\mathbb N}$ because $e^{\tau}\ge
\tau^m/m!$ for all $m\in{\mathbb N}$.

Note that $\psi(t)={\alpha t^{\frac{1}{n-1}}}$, $\alpha\in
(0,\infty)$, satisfies (\ref{eq3!!!33}) and $e^{\alpha
t^{\frac{1}{n-1}}}$ is convex on the segment $[T, \infty ]$ with
$T\ge\left[ \frac{n-2}{\alpha} \right]^{n-1}$, see e.g. the point
1.4.4 in \cite{Bou}.
\end{remark}

\medskip

\begin{corollary}\label{cor6.6.7A}{\sl\,
In particular, the conclusion of Theorem \ref{th4.2} is valid if,
for some $\alpha > 0$,
\begin{equation}\label{eq2!!!2}
\int\limits_D e^{\alpha Q^{\frac{1}{n-1}}(x)
}\frac{dm(x)}{\left(1+|x|^2\right)^n}\ \le\ M<\infty \ .
\end{equation}
The same is true for any function $\Phi=e^{\psi}$ where $\psi(t)$ is
a finite product of the function $\alpha t^{\beta}$, $\alpha >0$,
$\beta \ge 1/(n-1)$, and some of the functions
$[\log(A_1+t)]^{\alpha_1}$, $[\log\log(A_2+t)]^{\alpha_2},\ \ldots \
$, $\alpha_m\ge -1$, $A_m\in{\mathbb R }$, $m\in{\mathbb N}$, $t\in
[T,\infty ]$, $\psi(t)\equiv\psi(T)$,  $t\in [1,T]$ for a large
enough $T\in{\mathrm I}$.}
\end{corollary}

\medskip

\begin{remark}\label{remark4.3}
Note also that, as it follows from Lemma \ref{lem7.2}, the
conclusion of these theorems on con\-ver\-gen\-ce of ring
$Q$-homeomorpisms $f_j$ is also valid for ring
$Q_j$-ho\-meo\-mor\-pisms $f_j$ if at least one of the conditions on
$Q=Q_j$ in Theorems \ref{th6.6.1A}--\ref{th4.2} and Corollary
\ref{cor6.6.2A}--\ref{cor6.6.7A} holds uniformly in the parameter
$j=1,2,\ldots $. Thus, replacing the above conditions on $Q$ by the
corresponding uniform conditions on $K^{n-1}_{f_j}$, we obtain more
general results on convergence of the Sobolev homeomorphisms with no
point-wise dominant $Q(x)$ for $K^{n-1}_{f_j}(x)$. We give here some
of such results in the explicit form.\end{remark}

\bigskip


\begin{theorem}{}\label{th4.2.v}{\sl\,Let $D$ be a domain in ${\Bbb
R}^n,$ $n\ge 3,$ and let $\Phi:\overline{{\mathbb
R}^+}\rightarrow\overline{{\mathbb R}^+}$ be a nondecreasing convex
function satisfying the condition (\ref{eq3!}). Suppose that $f_j,$
$j=1,2,\ldots,$ is a sequence of ho\-me\-o\-mor\-phisms of $D$ into
${\Bbb R}^n$  in the class $W^{1,\varphi}_{\mathrm loc}$ with
(\ref{eqOS3.0a}) and
\begin{equation}\label{eq2!!.v}
\int\limits_D\Phi(K_{f_j}^{n-1}(x))\
\frac{dm(x)}{\left(1+|x|^2\right)^n}\ \le\ M<\infty\qquad\forall\
j=1,2,\ldots
\end{equation}
If $f_j\to f$ locally uniformly, then $f$ is either a constant in
$\overline{{\Bbb R}^n}$ or a ho\-meo\-mor\-phism into ${\Bbb R}^n.$}
\end{theorem}

\medskip

The notes in Remark \ref{rmkCOM3.1.1} to Theorem \ref{th4.2} on the
convexity of the function $\Phi$ and on its minorants are valid for
Theorem \ref{th4.2.v}. Thus, we have the following.

\medskip

\begin{corollary}\label{cor6.6.7Z}{\sl\,
In particular, the conclusion of Theorem \ref{th4.2.v} is valid if,
for some $\alpha > 0$,
\begin{equation}\label{eq2!!!2z}
\int\limits_D e^{\alpha K_f(x)
}\frac{dm(x)}{\left(1+|x|^2\right)^n}\ \le\ M<\infty \ .
\end{equation}
The same is true for any function $\Phi=e^{\psi}$ where $\psi(t)$ is
a finite product of the function $\alpha t^{\beta}$, $\alpha >0$,
$\beta \ge 1/(n-1)$, and some of the functions
$[\log(A_1+t)]^{\alpha_1}$, $[\log\log(A_2+t)]^{\alpha_2},\ \ldots \
$, $\alpha_m\ge -1$, $A_m\in{\mathbb R }$, $m\in{\mathbb N}$, $t\in
[T,\infty ]$, $\psi(t)\equiv\psi(T)$,  $t\in [1,T]$, with a large
enough $T\in{\mathrm I}$.}
\end{corollary}

\medskip

\begin{remark}\label{remark4.1}
Note that he integral conditions (\ref{eq2!!.v}) and (\ref{eq2!!})
for $K_f$ and $\Phi$ can be written in other forms that are
sometimes more convenient. Namely, by (\ref{eq1!!!}) with
$h(t)=t^{\frac{1}{n-1}}$ and $\varphi(t)=\Phi(t^{n-1})$, $\Phi =
\varphi\circ h$, the couple of the conditions (\ref{eq2!!.v}) and
(\ref{eq2!!}) is equivalent to the following couple
\begin{equation}\label{eq2!!!}
\int\limits_D\varphi (K_{f_j}(x))\
\frac{dm(x)}{\left(1+|x|^2\right)^n}\ \le\ M<\infty\qquad\forall\
j=1,2,\ldots
\end{equation}
and
\begin{equation}\label{eq3!!!}
\int\limits_{\delta}^{\infty} \frac{d\tau}{\tau\varphi^{-1}(\tau)}\
=\ \infty
\end{equation}
for some $\delta>\varphi(1).$ Moreover, by Theorem 2.1 in
\cite{RSY$_1$} the couple of the conditions (\ref{eq2!!!}) and
(\ref{eq3!!!}) for $\varphi = e^{\psi}$ is in turn equivalent to the
next couple
\begin{equation}\label{eq2!!!2}
\int\limits_D e^{\psi ( K_{f_j}(x))
}\frac{dm(x)}{\left(1+|x|^2\right)^n}\ \le\ M<\infty\qquad\forall\
j=1,2,\ldots
\end{equation}
and
\begin{equation}\label{eq3!!!3}
\int\limits_{\Delta}^{\infty}\ \psi(t)\ \frac{dt}{t^2}\ =\ \infty
\end{equation}
for some $\Delta>t_0$ where $t_0:=\sup\limits_{\psi(t)=-\infty}t,\ $
$t_0=1\ $ if $\ \psi(1)>-\infty .$ \end{remark}

\medskip

\cc
\section{ On normal families of homeomorphisms}  \label{16}

Here we give some criteria of normality for the Orlicz--Sobolev
classes of ho\-meo\-mor\-phisms.

\bigskip

Let $D$ be a domain in ${\Bbb R}^n$, $n\geqslant3$, and let
$\varphi:\overline{{\Bbb R}^+}\to\overline{{\Bbb R}^+}$ be a
nonconstant continuous nondecreasing function, $Q:D\to\bar{\mathrm
I}$ be a measurable function. Denote by ${\mathcal
O}_{Q,\Delta}^{\varphi}$ the collection of all homeomorphisms $f$ in
the Orlicz--Sobolev class $W^{1,\varphi}_{\mathrm loc}$ such that
$K^{n-1}_f(x)\leqslant Q(x)$ a.e. and $h\left(\overline{{\Bbb
R}^n}\setminus f(D)\right)$ $\geqslant\Delta>0$.  Moreover, denote
by ${\cal S}^p_{Q,\Delta}$ the class ${\mathcal
O}_{Q,\Delta}^{\varphi}$ with $\varphi(t)=t^p$, $p\geqslant1$.
Finally, let ${\mathcal S}_{Q ,\Delta}$ be the collection of all
homeomorphisms $f$ of $D$ into ${\Bbb R}^n$ in the class
$W^{1,1}_{\mathrm loc}$ such that $K^{n-1}_f(x)\leqslant Q(x)$ a.e.
and $h\left(\overline{{\Bbb R}^n}\setminus
f(D)\right)\geqslant\Delta>0$. Here and further $h(A)$ denote the
spherical (chordal) diameter of a set $A$ in $\overline{{\Bbb
R}^n}$.

\bigskip

\begin{lemma}\label{lem7.3}
{\sl\, Let $D$ be a domain in ${\Bbb R}^n,$ $n\ge 3,$ $\Delta>0$ and
$Q: D\rightarrow \bar{\mathrm I}$ be a measurable function with the
condition (\ref{eq6.3.24}). Then the classes ${\cal
O}_{Q,\Delta}^{\varphi}$, ${\cal S}^p_{Q,\Delta}$ and ${\cal
S}_{Q,\Delta}$ form normal families if $\varphi$ satisfies
(\ref{eqOS3.0a}), correspondingly, $p>n-1$ and $Q\in
L^{\gamma}_{\mathrm loc}$, $\gamma >1$. }
\end{lemma}

\bigskip

\begin{remark}\label{remark6.1} Lemma \ref{lem7.3} follows directly
from Lemma 6.1 in \cite{RS$_3$} for $Q$-homeomorphisms. In
par\-ti\-cu\-lar, the conclusions of Lemma \ref{lem7.3} hold under
the conditions on $Q$ in Theorems \ref{th6.6.1A}--\ref{th4.2} and
Corollary \ref{cor6.6.2A}--\ref{cor6.6.7A}, cf. the corresponding
criteria of normality in \cite{KRSS}. See also Section 4 in
\cite{RS$_3$} for reasons on in\-ter\-con\-nec\-tions of various
conditions on $Q$ with the condition (\ref{eq6.3.24}).

Furthermore, as it follows from Remark 6.1 in \cite{RS$_3$}, the
conclusions on equi\-con\-ti\-nui\-ty and normality of ring
$Q$-homeomorpisms are valid for variable $Q$ but if at least one of
the conditions on $Q$ in Lemma \ref{lem7.2}, Theorems
\ref{th6.6.1A}--\ref{th4.2} and Corollary
\ref{cor6.6.2A}--\ref{cor6.6.7A} holds uniformly. Thus, again
replacing the above conditions on $Q$ by the corresponding uniform
conditions on $K^{n-1}_{f}$, we obtain more general results on
normality for classes of Sobolev's homeomorphisms with no point-wise
dominant $Q$ for $K^{n-1}_{f}$. We give here a few such results in
the explicit form.
\end{remark}

\bigskip

Now, let $\Phi:\bar{\mathrm I}\rightarrow\overline{{\mathbb R}^+}$
be a nondecreasing convex function. Given $\Delta\in (0,1)$ and
$M\in(0,\infty)$, denote by ${\mathcal O}_{\Phi}^{\varphi}$ the
collection of all homeomorphisms $f$ of $D$ into ${\Bbb R}^n$ in the
Orlicz--Sobolev class $W^{1,\varphi}_{\mathrm loc}$ such that
$h\left(\overline{{\Bbb R}^n}\setminus f(D)\right)$
$\geqslant\Delta$ and
\begin{equation}\label{eq2!!.z}
\int\limits_D\Phi(K_{f}^{n-1}(x))\
\frac{dm(x)}{\left(1+|x|^2\right)^n}\ \le\ M\ .
\end{equation}
Further, let ${\cal S}^p_{\Phi}$ be the class ${\mathcal
O}_{\Phi}^{\varphi}$ with $\varphi(t)=t^p$, $p\geqslant1$.

Moreover, let ${\mathcal S}_{\Phi , \alpha}^*$, $\alpha\ge 1$, be
the collection of all homeomorphisms $f$ of $D$ into ${\Bbb R}^n$ in
the class $W^{1,1}_{\mathrm loc}$ such that $h\left(\overline{{\Bbb
R}^n}\setminus f(D)\right)\geqslant\Delta$ and
\begin{equation}\label{eq2!!.x}
\int\limits_D\Phi(K_{f}^{\alpha}(x))\
\frac{dm(x)}{\left(1+|x|^2\right)^n}\ \le\ M\ .
\end{equation}

Finally, given a nondecreasing function $\psi:\bar{\mathrm
I}\to\overline{{\Bbb R}^+}$ such that ${\psi^m\left(
t^{\frac{1}{n-1}}\right)}$ is convex for some $m\in{\mathbb N}$ on a
segment $[T,\infty ]$ for some $T\in{\mathrm I}$, denote by
${\mathcal S}^{\psi}_*$ the collection of all homeomorphisms $f$ of
$D$ into ${\Bbb R}^n$ in the class $W^{1,1}_{\mathrm loc}$  such
that $h\left(\overline{{\Bbb R}^n}\setminus
f(D)\right)\geqslant\Delta$ and
\begin{equation}\label{eq2!!.y}
\int\limits_De^{\psi(K_{f}(x))}\
\frac{dm(x)}{\left(1+|x|^2\right)^n}\ \le\ M\ .
\end{equation}
Apply to the class ${\mathcal S}^{\psi}_*$ Remark 6.1 in
\cite{RS$_3$} side by side with Remark \ref{rmkCOM3.1.1} above.

\medskip

\begin{lemma}\label{lem7.3.z}
{\sl\, Let $D$ be a domain in ${\Bbb R}^n,$ $n\ge 3$. Then the
classes ${\cal O}_{\Phi}^{\varphi}$ and ${\mathcal S}^{\psi}_*$,
${\cal S}^p_{\Phi}$ and ${\cal S}_{\Phi , \alpha}^*$ form normal
families if $\varphi$, $\Phi$ and $\psi$ satisfies (\ref{eqOS3.0a}),
(\ref{eq3!}) and (\ref{eq3!!!3}), correspondingly, $p>n-1$ and
\begin{equation}\label{eq7.77.1.j} \int\limits_{\delta}^{\infty}
\frac{d\tau}{\tau\left[\Phi^{-1}(\tau)\right]^{\frac{1}{\alpha}}}\
=\ \infty\ , \ \ \ \ \ \delta>\Phi(1)\ ,\ \ \  \alpha >n-1\ .
\end{equation}
}
\end{lemma}

\medskip

\cc
\section{ Main Lemma }  \label{16}

The following lemma was proved earlier for mappings with bounded
distortion in the paper \cite{GMRV}, Lemma 4.7, and for mappings
with finite length distortion in the monograph \cite{MRSY}, Lemma
8.6.

\medskip

\begin{lemma}
{}\label{le7.8.2}{\sl\, Let $D$ be a domain in ${\Bbb
R}^{n},$\thinspace \thinspace $n\ge 2,$ and let $f_{j}$
$j=1,2,\ldots ,$ be a sequence of ho\-me\-o\-mor\-phisms of $D$ into
${\Bbb R}^{n}$ of the class $W^{1,1}_{\mathrm loc}$ converging
locally uniformly to a mapping $f:D\rightarrow {\Bbb R}^{n}.$ Then
\begin{equation}
P_{f}(x_{0})\le\liminf\limits_{h\rightarrow 0}\text{ }%
\liminf\limits_{j\rightarrow \infty }\,\frac{1}{h^{n}}\int%
\limits_{C(x_{0},h)}P_{f_j}(y)\ dm(y)  \label{eq7.8.1A}
\end{equation}%
at each point $x_{0}$ of differentiability of the mapping $f$ where
$C(x_{0},h)$ denotes the cube in ${{\Bbb R}^{n}}$ centered at
$x_{0}$
whose edges of the length $h$ are oriented along the principal axes of the quadratic form $%
(f^{\,\prime }(x_{0})z,f^{\,\prime }(x_{0})z)$.}
\end{lemma}

\medskip

Here we use the dilatation $P_f$ defined in (\ref{eqOS1.3.1}).

\medskip

\begin{proof} Without los of generality we may assume that $x_0=0,$ $f(0)=0,$
and that $f_{j}(0)=0$, $j=1,2,\ldots $ . Let $e_{1},\ldots ,e_{n}$
be an
orthonormal basis in ${{\Bbb R}^{n}}$ formed by the eigenvectors of $%
f^{\prime }(0)^{\ast }f^{\prime }(0)$. Note that $f^{\prime
}(0)\mathbb{B}^{n}$ is an ellipsoid whose semiaxes $\lambda _{1}\le
\lambda _{2}\le \dots \le \lambda _{n}$ are the nonnegative square
roots of the corresponding eigenvalues of $f^{\prime }(0)^{\ast
}f^{\prime }(0)$. We may assume that $f^{\prime }(0)\ne 0$, i.e.,
$\lambda_n>0,$ because in the contrary case $P_f(0)=1$ and the
inequality (\ref{eq7.8.1A}) is obvious. We also abbreviate
$C(h)=C(0;h)$.

Now, for every $\varepsilon
>0,$ we can choose $\delta =\delta (\varepsilon )>0$ such that for
$h\in (0,\delta )$ and all $y\in C(h),$
$$
|f(y)-f^{\prime }(0)y|<h\varepsilon
$$%
because $f$ is differentiable at $0$. Moreover, since
$f_{j}\rightarrow f$ locally uniformly, we have, for all $y\in
C(h),$ that
\begin{equation}
|f_{j}(y)-f^{\prime }(0)y|<h\varepsilon  \label{eq7.8.9a}
\end{equation}%
under $j>j_{0}$. The set $f^{\prime }(0)C(h)$ is the rectangular
parallelopiped
$$
(-\lambda _{1}h/2,\lambda _{1}h/2)\times \cdots \times (-\lambda
_{n}h/2,\lambda _{n}h/2),
$$%
that can be degenerate if $\lambda_i=0$, $i=1,\ldots , k$, $0<k<n$,
whose edges are oriented along the basis vectors $\tilde{e}_{k+1},\ldots ,%
\tilde{e}_{n}$, of ${{\Bbb R}^{n}}$,
$$
\tilde{e}_{i}={\frac{f^{\prime }(0)e_{i}}{|f^{\prime
}(0)e_{i}|}}\,,\text{ \ \ }\ i=k+1,\ldots ,n,
$$%
and some vectors $\tilde{e}_{1},\ldots ,%
\tilde{e}_{k}$ that form orthonormal basis  of the orthogonal
complement to $(n-k)$-dimensional subspace of ${{\Bbb R}^{n}}$
generated by $\tilde{e}_{k+1},\ldots ,%
\tilde{e}_{n}$. Inequality (\ref{eq7.8.9a}) yields that all the
points $f_{j}(y)$, $y\in C(h)$, lie in the parallelopiped
$$
\left( -\left( {\frac{\lambda _{1}}{2}}+\varepsilon \right) h,\ \left( {%
\frac{\lambda _{1}}{2}}+\varepsilon \right) h\right) \times \cdots
\times
\left( -\left( {\frac{\lambda _{n}}{2}}+\varepsilon \right) h,\ \left( {%
\frac{\lambda _{n}}{2}}+\varepsilon \right) h\right) .
$$%
Here ${{\Bbb R}^{n}}$ is again equipped with the basis $\tilde{e}%
_{1},\ldots ,\tilde{e}_{n}$. Thus,
\begin{equation}
\text{mes}(f_{j}(C(h)))\le h^{n}(\lambda _{1}+2\varepsilon )(\lambda
_{2}+2\varepsilon )\cdots (\lambda _{n}+2\varepsilon ).
\label{eq7.8.10a}
\end{equation}
and by (\ref{eqOS1.1.1}) we have the inequalities
\begin{equation}
\int\limits_{C(h)}|J_{f_{j}}(y)|\ dm(y)\ \le \
\text{mes}(f_{j}(C(h)))\le h^{n}(|J_f(0)|+\triangle (\varepsilon )),
\label{eq7.8.11a}
\end{equation}%
where $\triangle (\varepsilon )\rightarrow 0$ as $\varepsilon
\rightarrow 0$ because $J_f(0)=\lambda _{1}\lambda _{2}\cdot \ldots
\cdot \lambda _{n}$.
\medskip

Next, consider the $(n-1)$-dimensional cube $C^{\ast }(h)$ with center at $%
x=0 $ and edges of the length $h$ oriented along $e_{1},\ldots
,e_{n-1}$. Consider a segment $l(z)$, $z\in C^{\ast }(h)$,
perpendicular to $C^{\ast
}(h)$ inside of $C(h)$ and write $l_{j}(z)$ for the length of the path $%
f_{j}(l(z))$. Since $f_j\in W_{\rm loc}^{1,1},$ by Theorem 1 of
Section 1.1.3 and Theorem of Section 1.1.7 in \cite{Maz} $f_j$ is
absolutely continuous on $l(z)$ for a.e. $z\in C^{\ast }(h)$ and by
Theorem 1.3 in \cite{Va} we have that
\begin{equation}
l_{j}(z)\ =\ \int\limits_{-h/2}^{h/2}|f_{j}^{\prime
}(z,y_{n})e_{n}|\ dy_{n} \label{eq7.8.12a}
\end{equation}%
for almost every $z\in C^{\ast }(h)$ with respect to the
$(n-1)$-dimensional Lebesgue measure. On the other hand,
(\ref{eq7.8.9a}) implies that
$$
l_{j}(z)\ge (|f^{\prime }(0)e_{n}|-2\varepsilon )h=(\lambda
_{n}-2\varepsilon )h.
$$%
Hence (\ref{eq7.8.12a}) yields that
$$
\int\limits_{-h/2}^{h/2}|f_{j}^{\prime }(z,y_{n})e_{n}|\ dy_{n}\ \ge
\ h(\lambda _{n}-2\varepsilon )
$$%
for a.e. $z\in C^{\ast }(h)$. Thus, integrating over $C^{\ast }(h)$
and using the Fubini theorem, we obtain that
\begin{equation}
\int\limits_{C(h)}|f_{j}^{\prime }(y)e_{n}|\ dm(y)\ \ge \
h^{n}(\lambda _{n}-2\varepsilon ).  \label{eq7.8.13a}
\end{equation}

Let first $K_{f_j}(y)\ne \infty$ a.e. in $C(h).$ Then the H\"{o}lder
inequality gives
\begin{equation}\label{eq7.8.14a}
\int\limits_{C(h)}|f_{j}^{\prime }(y)e_{n}|\ dm(y) \le
\int\limits_{C(h)}\Vert f_{j}^{\prime }(y)\Vert \ dm(y)\ =\
\int\limits_{C(h)}K_{f_{j}}^{1/n}(y)\ J_{f_{j}}(y)^{1/n}\ dm(y)\
\end{equation}
$$\le \left( \int\limits_{C(h)}K_{f_{j}}^{{\frac{1}{{n-1}}}}(y)\ dm(y)\right) ^{
\frac{{n-1}}{n}}\left( \int\limits_{C(h)}J_{f_{j}}(y)\ dm(y)\right)
^{{\frac{1}{n} }}.$$ Here the equality $\Vert f_{j}^{\prime
}(y)\Vert ^{n}=K_{f_{j}}(y)\ J_{f_{j}}(y)$
a.e. has also been used. Combining (\ref{eq7.8.14a}), (\ref{eq7.8.11a}%
) and (\ref{eq7.8.13a}), we obtain that
$$
\left( {\frac{{(\lambda _{n}-2\varepsilon
)^{n}}}{{|J_f(0)|+\triangle
(\varepsilon )}}}\right) ^{{\frac{1}{{n-1}}}}\le {\frac{1}{h^{n}}}%
\int\limits_{C(h)}K_{f_{j}}^{{\frac{1}{{n-1}}}}(y)\ dm(y)\ .
$$%

Note that the last inequality also  holds in the evident way for the
case $K_{f_j}(y)= \infty$ on a set of positive measure because in
this case the right hand side is equal to
$\infty$.

\bigskip

Thus, letting here first $j\rightarrow \infty $, then $h\rightarrow 0$ and, finally, $%
\varepsilon \rightarrow 0,$ we complete the proof. \end{proof}
$\Box$

\bigskip

Applying the Jensen inequality to (\ref{eq7.8.1A}), we obtain the
following conclusion.

 \medskip

\begin{corollary}
{}\label{cor7.8.1}{\sl Under the hypotheses and notations of Lemma
\ref{le7.8.2},
\begin{equation}
\,\,\,\,\,\,\,\,\,\,\,\,\,\,\Phi \left( P_{f}\left( x_{0}\right)
\right)
\,\,\le \,\,\liminf\limits_{h\rightarrow 0}\text{ }\liminf\limits_{j%
\rightarrow \infty }\frac{1}{h^{n}}\int\limits_{C\left(
x_{0},h\right) }\Phi \left( P_{f_j}\left( y \right) \right) \ dm(y)
\label{eq7.8.12}
\end{equation}%
for every continuous convex function $\Phi :\bar{\mathrm I}
 \to \overline{\mathbb{R}^+}$.}
\end{corollary}

 \medskip

In particular, for the function $\Phi(t)=t^{n-1}$, we have the next
conclusion.

 \medskip

\begin{corollary}
{}\label{cor7.8.2}{\sl Under the hypotheses of Lemma \ref{le7.8.2},
\begin{equation}
K_{f}(x_{0})\,\,\le \,\,\liminf\limits_{h\rightarrow 0}\text{ }%
\liminf\limits_{j\rightarrow \infty }\,\frac{1}{h^{n}}\int%
\limits_{C(x_{0},h)}K_{f_{j}}(y)\,\,dm(y)\,.  \label{eq7.8.13}
\end{equation}}
\end{corollary}

 \medskip

 \begin{theorem}
{}\label{th7.8.2}{\sl Let $D$ be a domain in ${\Bbb
R}^{n},$\thinspace \thinspace $n\ge 2,$ and let $f_{j}$,
$j=1,2,\ldots ,$ be a sequence of ho\-me\-o\-mor\-phisms of $D$ into
${\Bbb R}^{n}$ of the class $W^{1,1}_{\mathrm loc}$ converging
locally uniformly to a mapping $f:D\rightarrow {\Bbb R}^{n}$ which
is differentiable
 a.e. in $D.$ If
\begin{equation}
P_{f_j}\left(x\right) \,\,\le \,\,P\left(x\right)\,\in L_{\mathrm{loc%
}}^{1}\ \quad j=1,2,\ldots ,  \label{eq7.8.2}
\end{equation}%
then
\begin{equation}
P_{f}\left(x\right) \,\le \,\limsup\limits_{j\rightarrow \infty
}\,P_{f_j}\left(x\right) \qquad{\rm a.e.}  \label{eq7.8.3}
\end{equation}}
Moreover, if $f$ is a homeomorphism, then $f\in W^{1,1}_{\mathrm
loc}$ and $\partial_if_j\to\partial_if$ weakly in $L^1_{\mathrm
loc}$ as $j\to\infty$ for all $i=1,\ldots n$.
\end{theorem}

\medskip

\begin{proof} Applying Lemma \ref{le7.8.2}, the condition (\ref{eq7.8.2}) and the theorem on
term-by-term integration, see e.g. Theorem I.12.12 in \cite{Sa}, we
have that
\begin{equation}
P_{f}\left(x\right) \,\,\le \,\liminf\limits_{h\rightarrow 0}\,\frac{1}{%
h^{n}}\int\limits_{C\left( x,h\right) }\limsup\limits_{j\rightarrow
\infty }P_{f_j}(y)\,\,dm(y).  \label{eq7.8.4}
\end{equation}%
Now, by the  theorem on the differentiability of the indefinite
Lebesgue integral, see e.g. Theorem IV.6.3 in \cite{Sa}, we obtain
a.e. the equality
\begin{equation}
\qquad\lim\limits_{h\rightarrow 0}\,\frac{1}{h^{n}}%
\int\limits_{C\left( x,h\right) }\limsup\limits_{j\rightarrow \infty
}P_{f_j}\left(y\right) \,\,dm(y)\,=\limsup\limits_{j\rightarrow
\infty }P_{f_j}\left(x\right) .\,\,  \label{eq7.8.5}
\end{equation}%
Finally, combining $(\ref{eq7.8.4})$ and $(\ref{eq7.8.5}),$ we come to $(\ref{eq7.8.3}%
).$

Next, by (\ref{eq7.8.2}) $f_j$ are of finite distortion,
$P_{f_j}<\infty$ a.e., and by (\ref{eqOSadm.1}) , for every compact
set $C\subset \Omega$, $i=1,\ldots , n$, $j=1,2,\ldots$
\begin{equation}\label{eqOSadm.1.1.a}
\Vert\partial_i f_j\Vert\ \leqslant\ \Vert
P_{f_j}\Vert^{(n-1)/n}\cdot|f_j(C)|^{1/n}\ \leqslant\ \Vert
P\Vert^{(n-1)/n}\cdot|f_j(C)|^{1/n}\  <\ \infty
\end{equation}
where $\Vert * \Vert$ denotes the norm of $*$ in $L^1(C)$. Take
$\rho_*\in(0,\rho)$ where $\rho=\mbox{dist}(C,\partial D)$ and cover
$C$ by all balls $\mathrm B(x,\rho_*)$, $x\in C$. Choosing a finite
covering from the given covering, we obtain an open set $V$
containing $C$ and such that $\overline{V}\subset D$. By the
construction, $\overline{V}$ is a compact subset of $D$ and,
moreover, $f(C)$ and $f(\overline{V})$ are compact subsets of the
domain $D^{\prime}=f(D)$ such that $f(C)\subset f(V)$, $f(V)$ is
open and hence $\mbox{dist}(f(C),
\partial f(V))>0$ because $f$ is a homeomorphism. Consequently,
$f_j(C)\subset f(V)$ for large enough $j$ and we have by
(\ref{eqOSadm.1.1.a}) that
\begin{equation}\label{eqOSadm.1.1.1.a}
\Vert\partial_i f_j\Vert\  \leqslant\ \Vert
P\Vert^{(n-1)/n}\cdot|f(\overline{V})|^{1/n}\ <\ \infty
\end{equation}
for such $j$. It is proved similarly that
\begin{equation}\label{eqOSadm.1.1.1.b}
\Vert\partial_i f_j\Vert_E\  \leqslant\ \Vert
P\Vert^{(n-1)/n}_E\cdot|f(\overline{V})|^{1/n}\ <\ \infty
\end{equation}
where $\Vert * \Vert_E$ denotes the norm of $*$ in $L^1(E)$ for an
arbitrary measurable set $E\subseteq C$. Thus, $f\in
W^{1,1}_{\mathrm loc}$ and $\partial_if_j\to\partial_if$ weakly in
$L^1_{\mathrm loc}$, $i=1,\ldots n$ as $j\to\infty$ by Lemma 2.1 in
\cite{RSY}. $\Box$ \end{proof}

\medskip

In the next section we give the corresponding results with no
point-wise dominant $P(x)$.

\cc
\section{ On Semicontinuity of Dilatations in the Mean}  \label{16.1}


The corresponding analogs of the following result for the plane case
can be found in the monograph \cite{GR^*}, see Theorem 11.1 therein,
and for space mappings with bounded distortion in the paper
\cite{GMRV}, Theorem 4.1.

\medskip

As above, we use here the dilatation $P_f$ defined in
(\ref{eqOS1.3.1}).

\medskip

\begin{theorem}
{}\label{th7.8.2.1}{\sl Let $D$ be a domain in ${\Bbb
R}^{n},$\thinspace \thinspace $n\ge 2,$ and let $f_{j}$,
$j=1,2,\ldots ,$ be a sequence of ho\-me\-o\-mor\-phisms of $D$ into
${\Bbb R}^{n}$ of the class $W^{1,1}_{\mathrm loc}(D)$ converging
locally uniformly to a mapping $f:D\rightarrow {\Bbb R}^{n}$ which
is differentiable
 a.e. in $D.$ Then on every open set $\Omega\subseteq D$
\begin{equation}\label{eq_1.14}
    \int\limits_{\Omega} \Phi (P_f (x))\ dm(x)\ \leq\
    \liminf\limits_{j \to \infty} \int\limits_{\Omega} \Phi (P_{f_j} (x))\
    dm(x)
\end{equation}
for every strictly convex function $\Phi : \bar{\mathrm I} \to
\overline{\mathbb{R}^+}$ which is continuous from the left at the
point $T = \sup\limits_{\Phi(t)<\infty}\ t$. Moreover, if $f$ is a
homeomorphism and the right hand side in (\ref{eq_1.14}) is finite,
then $f\in W^{1,1}_{\mathrm loc}(\Omega)$. Finally, if in addition
\begin{equation}\label{eq_1.14.z}
    \limsup\limits_{j \to \infty} \int\limits_{\Omega} \Phi (P_{f_j} (x))\
    dm(x)\ <\ \infty\ ,
\end{equation}} then $\partial_if_{j}\to\partial_if$ weakly in
$L^1_{\mathrm loc}(\Omega)$ as $j\to\infty$ for all $i=1,\ldots
n$.\end{theorem}

\medskip

\begin{remark}\label{remark5.1} Here the continuity of the function $\Phi$
is understood in the sense of the topology of
$\overline{\mathbb{R}^+}$. The conditions that $\Phi$ is convex,
nondecreasing and continuous from the left at the point $T$ are not
only sufficient but also necessary to conclude (\ref{eq_1.14}) as it
was already shown for the plane case in the monograph \cite{GR^*}.
We will prove this fact for the space case in the final section.
\end{remark}

\medskip

Before the proof of Theorem \ref{th7.8.2.1}, let us give a
reformulation, in terms of series, of the known Fatou lemma, see
e.g. Theorem I.12.10 in \cite{Sa}.

\medskip

\begin{lemma}
{}\label{le7.8.2.1}{\sl\, Let $a_{m j}, \ m, \, j = 1, \, 2, \,
\dots $ be nonnegative numbers. Then
\begin{equation}\label{eq_1.29}
    \sum\limits_{m=1}^{\infty}\ \liminf\limits_{j \to \infty}\ a_{m
    j}\
    \leq\ \liminf\limits_{j \to \infty}\ \sum\limits_{m=1}^{\infty}\
    a_{m j}\ .
\end{equation}
If a sequence $a_{m j}$, $ m, \, j = 1, \, 2, \, \dots $ contains
negative numbers, then to have (\ref{eq_1.29}) it is sufficient in
addition to assume that $|a_{m j}| \leq b_m$ where $\sum b_m <
\infty$.}
\end{lemma}

\medskip

Indeed, the case of series is reduced to the standard integral case
through applying the step--functions $\varphi_j$, $j= 1, \, 2, \,
\dots$ given on the segment $[0, 1]$\ :
$$
\varphi_j(t)\ =\ 2^m a_{m j}\ ,\ \ \ \ \ \  \sum\limits_{k=1}^{m-1}\
2^{-k}\ \leq\ t\ <\ \sum\limits_{k=1}^{m}\ 2^{-k}, \quad m=1, \, 2,
\, \dots \
$$

\medskip

\begin{proof} If the right hand side in (\ref{eq_1.14}) is equal to
$\infty$, then the inequality (\ref{eq_1.14}) is evident. Hence we
may assume further that
\begin{equation}\label{eq_1.14.1}
   \int\limits_{\Omega} \Phi (P_{f_j} (x))\
    dm(x)\ \le\ M\ <\ \infty\ \ \ \ \ \ \ \  \forall\ j=1,2,\ldots\ .
\end{equation}

1) Let first the function $\Phi (P_f(z))$ is locally integrable in a
bounded $\Omega$. Then by the  theorem on the differentiability of
the indefinite integral, IV(6.3) in \cite{Sa},
$$
\lim\limits_{h \to 0} \, \frac{1}{h^2} \int\limits_{C(x; h)} \Phi
(P_f (\zeta) \, dm(\zeta) = \Phi (P_f(x))
$$
and, moreover, (\ref{eq7.8.12}) holds by Corollary \ref{cor7.8.1}
for $x \in E$ where $|\Omega \setminus E| = 0.$ Thus,
$$
\int\limits_{C(x, h)} \Phi (P_f(\zeta)) \, dm(\zeta) \leq
\liminf\limits_{n \to \infty} \iint\limits_{C(x, h)} \Phi
(P_{f_j}(\zeta)) \, dm(\zeta) + \varepsilon h^2,
$$
for every point $x \in E$ and every $\varepsilon
> 0$ with small enough $h < \delta = \delta (\varepsilon, x)$.

The system of cubes $C (x, h),$  $x\in E,$  $h< \min (\rho (x,
\partial \Omega) / \sqrt{n}, \delta (\varepsilon, x))$, forms
a covering of the set $E$ in the sense of Vitali and we are able to
choose a countable collection of mutually disjoint cubes $C_m= C
(x_m, h_m)$ in the given covering such that $|E \setminus \cup E_m|
= 0$, $C_m \subseteq \Omega$ and $|\Omega \setminus \cup E_m| =0,$
$|\Omega| = \sum |E_m|$, see e.g. Theorem IV.3.1 in \cite{Sa}.

By countable additivity of the integral, see e.g. Theorem I.12.7 in
\cite{Sa}, and by Lemma \ref{le7.8.2.1} applied to
$$
a_{m j}= \int\limits_{E_m} \Phi (P_{f_j} (\zeta)) \, dm(\zeta) ,
\quad m, \, j= 1, \, 2, \, \dots\,,
$$
we obtain that
$$
\int\limits_{\Omega} \Phi (P_f (\zeta)) \, dm(\zeta) \leq
\liminf\limits_{j \to \infty} \int\limits_{\Omega} \Phi (P_j
(\zeta))\,  dm(\zeta) + \varepsilon |\Omega|
$$
and by arbitrariness of $\varepsilon > 0$ we come to
(\ref{eq_1.14}).

\bigskip

2) If $\Phi(T)<\infty$, then also $T<\infty$ because of the
condition $\Phi(t)/t\to\infty$ as $t\to\infty$. Consequently,
$P_{f_j}(x)\le T$ a.e. by the condition (\ref{eq_1.14.1}) and hence
 $P_{f}(x)\le T$ a.e. by Lemma \ref{le7.8.2}. However, then $\Phi(P_{f}(x))\le
 \Phi(T)$ a.e. and $\Phi(P_{f}(x))$ is locally integrable in
 $\Omega$. Thus, by the point 1) of the proof we obtain (\ref{eq_1.14}) for bounded $\Omega$.

\bigskip

3) Now, let $\Phi(T)=\infty$. Then by continuity of $\Phi$ in the
sense of $\overline{\mathbb{R}^+}$ from the left at the point $T$,
there is a monotone sequence $t_m\in(1, T)$, $\ m= 1, \, 2, \,
\dots$ $(1, \infty)$ such that $t_m \to T$, correspondingly,
$\Phi(t_m) \to \infty$ as $m \to \infty$ and $\Phi (t)$ is
differentiable at every point $t_m$, $\Phi' (t_m)
> 0$, see e.g. Corollary I.4.2 in \cite{Bou}. Let us consider the functions
$$
\Phi_m(t) = \left\{
  \begin{array}{lr}
    \Phi (t), &  t \leq t_m, \\
    \Phi (t_m) + \Phi' (t_m) (t- t_m),   & \quad t \geq t_m;
  \end{array}
\right.
$$
and
$$
\varphi_m(\tau)= \left\{
  \begin{array}{lr}
    \tau, & \tau \leq \Phi (t_m), \\
     \Phi (\alpha_m + \beta_m \tau),  & \quad \tau \geq \Phi(t_m);
  \end{array}
\right.
$$
where the coefficients
$$
\left\{
  \begin{array}{ll}
    \alpha_m= t_m - \Phi (t_m) / \Phi' (t_m) & \\
    \beta_m = 1 / \Phi'(t_m) &
   \end{array}
\right.
$$
have been found from the condition
\begin{equation}\label{eq_1.34.1}
\alpha_m + \beta_m \left[ \Phi(t_m) + \Phi'(t_m) (t- t_m) \right]
\equiv t
\end{equation}
and, thus, by the construction
\begin{equation}\label{eq_1.34}
    \varphi_m (\Phi_m(t)) \equiv \Phi (t), \quad m= 1, \, 2, \, \dots \,.
\end{equation}
As it is easy to see, all the functions $\Phi_m$ and $\varphi_m$ are
nondecreasing, continuous and convex, see e.g. Proposition I.4.8 in
\cite{Bou}, and
\begin{equation}\label{eq_1.35}
    \lim\limits_{\tau \to \infty} \frac{\varphi_m (\tau)}{\tau} =
    \infty
\end{equation}
because $\Phi$ is continuous and strictly convex. Moreover, the
sequence $\Phi_m\le \Phi$ is increasing and  point-wise convergent
to $\Phi$ as $m \to \infty$.

\medskip

Let us first prove that  the functions  $\Phi_m (P\,(x)), \ m=1, \,
2, \, \dots,$ are integrable over bounded $\Omega$. Indeed, by the
Jensen inequality, (\ref{eq_1.14.1}) and (\ref{eq_1.34}), for every
$E \subseteq \Omega$ with $|E| > 0$, we may consider that for all
$j=1, \, 2, \, \dots$ and $M_{*}= M + \varepsilon$:
\begin{equation}\label{eq_1.36}
    \varphi_m \left( \frac{1}{|E|} \int\limits_{E} \Phi_m (P_{f_j}(\zeta)) \,
    dm (\zeta)
    \right) \leq \frac{M_{*}}{|E|} < \infty.
\end{equation}

\medskip

If $T < \infty$, then by definition of $\varphi_m$,
(\ref{eq_1.34.1}) and (\ref{eq_1.36}) we obtain that
$$
\frac{1}{|E|} \int\limits_{E} \Phi_m (P_{f_j}(\zeta)) \ dm(\zeta)\
\leq\ Q_m\ \colon =\ \Phi (t_m) + \Phi'(t_m) (T - t_m)
$$
and hence $\Phi_m (P_{f_j}(x)) \leq Q_m $, $\ m, \, j=1, \, 2, \,
\dots, $ for a.e. $x \in \Omega$ by the theorem on differentiability
of the indefinite integral. Thus, by Theorem \ref{th7.8.2}  $\Phi_m
(P(x)) \leq Q_m$ a.e. and the lntegrability of $\Phi_m (P)$ become
evident.

\medskip

If $T = \infty$, then, in view of (\ref{eq_1.35}) and
(\ref{eq_1.36}),
$$
\int\limits_{E} \Phi_m (P_{f_j}(\zeta))\ dm(\zeta) \leq M_{*} \,
\frac{\tau} {\varphi_m(\tau)} \to 0,
$$
where $\tau = \varphi_m^{-1} (M_{*} / |E|) \to \infty$ as $ |E| \to
0$. Moreover,
$$
\| \Phi_m (P_{f_j}) \|_{L^1 (\Omega)} \leq |\Omega| \,
\varphi^{-1}_m (M_{*} / |\Omega|).
$$
Thus, the sequence $\Phi_m(P_{f_j})|_{\Omega}$ is weakly compact in
$L^1 (\Omega)$, see e.g. Corollary IV.8.11 \cite{DS}. Thus, we may
assume that $\Phi_m(P_{f_j})|_{\Omega}\to\Psi$ weakly in $L^1
(\Omega)$ as $j \to \infty$. By Corollary \ref{cor7.8.1} and the
theorem on the differentiability of the indefinite integral, see
e.g. Theorem IV.6.3 in \cite{Sa}, we obtain that
$\Phi_m(P_f)\le\Psi$ a.e. in $\Omega$ and, consequently,
$\Phi_m(P_f)|_{\Omega} \in L^1 (\Omega)$.

\medskip

Hence we may apply the first point of the proof to every of the
functions $\Phi_m,$ $m= 1, \, 2, \, \dots$ to obtain the inequality
$$
\int\limits_{\Omega} \Phi_m (P_f(\zeta)) \, dm(\zeta) \leq
\liminf\limits_{j \to \infty} \int\limits_{\Omega} \Phi_m
(P_{f_j}(\zeta)) \, dm(\zeta).
$$
By the Lebesgue theorem on integration of monotone sequences of
functions, see e.g. Theorem IV.12.6 in \cite{Sa},
$$
\lim\limits_{m \to \infty} \int\limits_{\Omega} \Phi_m (P_f(\zeta))
\ dm(\zeta)\ =\ \int\limits_{\Omega} \Phi (P_f(\zeta)) \ dm(\zeta)\
,
$$
$$
\lim\limits_{m \to \infty} \int\limits_{\Omega} \Phi_m
(P_{f_j}(\zeta))\ dm(\zeta)\ =\ \int\limits_{\Omega} \Phi
(P_{f_j}(\zeta)) \ dm(\zeta)\ .
$$
It remains to note that the double sequence of numbers
$$
a_{m j}= \int\limits_{\Omega} \Phi_m (P_{f_j}(\zeta))\  dm(\zeta)
$$
is increasing in  $m$ and hence, see Lemma \ref{le7.8.2.1},
\begin{equation}\label{eq_1.37}
    \lim\limits_{m \to \infty} \lim\limits_{\overline{j \to \infty}} a_{m j} \leq
    \lim\limits_{\overline{j \to \infty}} \lim\limits_{m \to \infty} a_{m j}
\end{equation}
that leads us to the relation (\ref{eq_1.14}) in the case of bounded
 $\Omega$.

\medskip

4) Finally, applying the exhausting $\Omega_m= \left\{ \zeta \in
\Omega : |\zeta | < m \right\},$ $m=1, \, 2, \, \dots$ and the
inequality (\ref{eq_1.37}) from Lemma \ref{le7.8.2.1} to the other
double sequence
\begin{equation}\label{eq_1.38}
a_{m j}\ \colon =\ \int\limits_{\Omega_m} \Phi (P_{f_j}(\zeta)) \
dm(\zeta)\ =\ \int\limits_{\Omega} \chi_m (\zeta)\ \Phi
(P_{f_j}(\zeta)) \ dm(\zeta)\ ,
\end{equation}
where $\chi_m$ is the characteristic functions of the sets
$\Omega_m$, by the mentioned Lebesgue theorem we obtain
(\ref{eq_1.14}) in the general case.

\bigskip

5) By Theorem 3.1.2 in \cite{Rud}, since the function $\Phi$ is
strictly convex, we obtain from the condition (\ref{eq_1.14.z}) that
\begin{equation}\label{eq_1.38.v}
\int\limits_{\Omega}\ P_{f_j}(x)\ dm(x)\ \le\ K\ <\ \infty\ \ \ \ \
\ \forall\ j\ >\ N \end{equation} and, for every measurable set
$E\subset\Omega$,
\begin{equation}\label{eq_1.38.w} \int\limits_{E}\
P_{f_j}(x)\ dm(x)\ \le\ K_E\ <\ \infty\ \ \ \ \ \ \forall\ j\ >\ N
\end{equation} where
\begin{equation}\label{eq_1.38.y}
K_E\to 0\ \ \ \ \  \mbox{as}\ \ \ \ \ \ |E|\to 0\ . \end{equation}
By (\ref{eqOSadm.1}), (\ref{eq_1.38.v}) and (\ref{eq_1.38.v}), for
every compact set $C\subset \Omega$, $i=1,\ldots , n$, $j>N$,
\begin{equation}\label{eqOSadm.1.1.v}
\Vert\partial_i f_j\Vert\ \leqslant\ \Vert
P_{f_j}\Vert^{(n-1)/n}\cdot|f_j(C)|^{1/n}\ \leqslant\
K^{(n-1)/n}\cdot|f_j(C)|^{1/n}\  <\ \infty
\end{equation}
and
\begin{equation}\label{eqOSadm.1.1.w} \Vert\partial_i f_j\Vert_E\
\leqslant\ \Vert P_{f_j}\Vert_E^{(n-1)/n}\cdot|f_j(C)|^{1/n}\
\leqslant\ K_E^{(n-1)/n}\cdot|f_j(C)|^{1/n}\  <\ \infty
\end{equation}
where $\Vert * \Vert$ and $\Vert * \Vert_E$ denote the norm of $*$
in $L^1(C)$ and $L^1(E)$, correspondingly. Arguing similarly to the
proof of Theorem \ref{th7.8.2} we show that there is a compact set
$C_*$ such that $f_j(C)\subset f(C_*)$ for large enough $j$ and,
consequently, we have by (\ref{eqOSadm.1.1.v}) and
(\ref{eqOSadm.1.1.w}) that, for all $i=1,\ldots , n$,
\begin{equation}\label{eqOSadm.1.1.1.v}
\Vert\partial_i f_j\Vert\  \leqslant\
K^{(n-1)/n}\cdot|f(C_*)|^{1/n}\ <\ \infty\ \ \ \ \ \ \forall\ j\ >\
N_*
\end{equation}
and
\begin{equation}\label{eqOSadm.1.1.1.w}
\Vert\partial_i f_j\Vert_E\  \leqslant\
K^{(n-1)/n}_E\cdot|f(C_*)|^{1/n}\ <\ \infty\ \ \ \ \ \ \forall\ j\
>\ N_*
\end{equation}
By Lemma 2.1 in \cite{RSY} the conditions (\ref{eqOSadm.1.1.1.v})
and (\ref{eqOSadm.1.1.1.w}) imply that $f\in W^{1,1}_{\mathrm
loc}(\Omega)$ and $\partial_if_j\to\partial_if$ weakly in
$L^1_{\mathrm loc}(\Omega)$ as $j\to\infty$  for all $i=1,\ldots n$
whenever (\ref{eq_1.14.z}) holds. The proof is complete. $\Box$
\end{proof}

\bigskip

\begin{corollary}
{}\label{cor7.8.2.1}{\sl Under the hypotheses of Theorem
\ref{th7.8.2.1},
\begin{equation}\label{eq_1.14.111}
    \int\limits_{\Omega} \Phi (K_f (x))\ dm(x)\ \leq\
    \liminf\limits_{j \to \infty} \int\limits_{\Omega} \Phi (K_{f_j} (x))\ dm(x)
\end{equation}
for every continuous nondecreasing convex function $\Phi :
\overline{\mathbb{R}^+} \to \overline{\mathbb{R}^+}$ and every open
set $\Omega\subseteq D$.}
\end{corollary}

\bigskip

Let us give one of the most important consequence of Theorem
\ref{th7.8.2.1}.

\bigskip

\begin{corollary}
{}\label{cor7.8.2.1.2}{\sl Let $D$ be a domain in ${\Bbb
R}^{n},$\thinspace \thinspace $n\ge 2,$ and let $f_{j}$,
$j=1,2,\ldots ,$ be a sequence of ho\-me\-o\-mor\-phisms of $D$ into
${\Bbb R}^{n}$ of the class $W^{1,1}_{\mathrm loc}(D)$ converging
locally uniformly to a homeomorphism $f$ of $D$ into ${\Bbb R}^{n}$.
Then, for every $\alpha >n-1$ and every open set $\Omega\subseteq
D$,
\begin{equation}\label{eq_1.14.111.1}
    \int\limits_{\Omega} \Phi (K^{\alpha}_f (x))\ dm(x)\ \leq\
    \liminf\limits_{j \to \infty} \int\limits_{\Omega} \Phi (K^{\alpha}_{f_j} (x))\ dm(x)
\end{equation}
for every continuous nondecreasing convex function $\Phi
:\bar{\mathrm I} \to \overline{\mathbb{R}^+}$. Moreover, if $\Phi$
is not constant and the right hand side of (\ref{eq_1.14.111.1}) is
finite, then $f\in W^{1,p}_{\mathrm loc}(\Omega)$ with $p=n\alpha
/(1+\alpha)>n-1$ and $f$ is differentiable a.e. in $\Omega$.
Finally, if in addition
\begin{equation}\label{eq_1.14.111.1.z}
    \limsup\limits_{j \to \infty} \int\limits_{\Omega} \Phi (K^{\alpha}_{f_j} (x))\
    dm(x)\ <\ \infty\ ,
\end{equation}
then $\partial_if_{j}\to\partial_if$ weakly in $L^p_{\mathrm
loc}(\Omega)$ as $j\to\infty$ for all $i=1,\ldots n$.}
\end{corollary}

\medskip

\begin{proof} If the right hand side in (\ref{eq_1.14.111.1}) is
equal to $\infty$, then the inequality (\ref{eq_1.14.111.1}) is
evident. Hence we may assume, with no loss of generality, that
\begin{equation}\label{eq_1.14.1.a}
   \int\limits_{\Omega} \Phi (K^{\alpha}_{f_j} (x))\
    dm(x)\ \le\ M\ <\ \infty\ \ \ \ \ \ \ \  \forall\ j=1,2,\ldots\ .
\end{equation}
By the same reasons, we may also assume that $\Phi$ is not constant.
Finally, as in the last proof, by Lemma \ref{le7.8.2.1} and
countable additivity of the integral, we may assume that $\Omega$ is
a bounded domain.

To apply Theorem \ref{th7.8.2.1} we have to prove that $f$ is
differentiable a.e. in $\Omega$. Further we assume that the function
$\Phi$ is extended from $\bar{\mathrm I}$ to
$\overline{\mathbb{R}^+}$ by the equality $\Phi(t)\equiv \Phi(1)$
for all $t\in[0,1)$. Such extended function $\Phi$ is continuous,
nondecreasing, not constant and convex, see e.g. Proposition I.4.8
in \cite{Bou}. Setting $t_0=\sup\limits_{\Phi(t)=\Phi(0)} t$ and
$T_0=\sup\limits_{\Phi(t)<\infty} t$, we see that $t_0<T_0$ because
$\Phi$ is continuous, nondecreasing and not constant and, taking
$t_*\in(t_0,T_0)$, we have that
$$
\frac{\Phi(t)-\Phi(0)}{t}\ \ge \frac{\Phi(t_*)-\Phi(0)}{t_*}\ >\ 0\
\ \ \ \ \ \ \ \ \ \  \forall\ t\in [\ t_*\ ,\ \infty\ )
$$
by convexity of $\Phi$, see e.g. Proposition I.4.5 in \cite{Bou},
i.e., $\Phi(t)\ge at$ for $t\ge t_*$ where
$a={[\Phi(t_*)-\Phi(0)]}/{t_*}>0$. Thus, by (\ref{eq_1.14.1.a}) we
obtain that
\begin{equation}\label{eq_1.14.1.b}
   \int\limits_{\Omega} K^{\alpha}_{f_j} (x)\
    dm(x)\ \le\ t_*|\Omega|\ +\ M/a\ <\ \infty\ \ \ \ \ \ \ \  \forall\ j=1,2,\ldots\ .
\end{equation}
Hence by Proposition \ref{prOS2.4} $f_j\in W^{1,p}_{\mathrm
loc}(\Omega)$ with $p=n\alpha /(1+\alpha)>n-1$ because $\alpha>n-1$
and, more precisely, by (\ref{eqOSadm.1}) and (\ref{eq_1.14.1.b}),
for every compact set $C\subset \Omega$, $i=1,\ldots , n$,
$j=1,2,\ldots$
\begin{equation}\label{eqOSadm.1.1}
\Vert\partial_i f_j\Vert_p\ \leqslant\ \Vert
K_{f_j}\Vert^{1/n}_{\alpha}\cdot|f_j(C)|^{1/n}\ \leqslant\
A\cdot|f_j(C)|^{1/n}\ <\ \infty
\end{equation}
where $A=(t_*|\Omega| + M/a)^{1/\alpha}$. Arguing again as in the
proof of Theorem \ref{th7.8.2} we show that there is a compact set
$C_*$ such that $f_j(C)\subset f(C_*)$ for large enough $j$ and,
consequently, we have by (\ref{eqOSadm.1.1}) that
\begin{equation}\label{eqOSadm.1.1.1}
\Vert\partial_i f_j\Vert_p\  \leqslant\ A\cdot|f(C_*)|^{1/n}\ <\
\infty
\end{equation}
for such $j$. Thus, $\partial_if_{j}\to\partial_if$ weakly in
$L^p_{\mathrm loc}(\Omega)$ as $j\to\infty$ for all $i=1,\ldots n$
and $f\in W^{1,p}_{\mathrm loc}(\Omega)$ with $p>n-1$, see e.g.
Lemma 3.5, Ch. III, $\S\,3.4$ of the monograph \cite{Re}, and $f$ is
differentiable a.e. in $\Omega$  by the Gehring--Lehto--Menshov
theorem for $n=2$ and by the V\"{a}is\"{a}l\"{a} theorem for $n\ge
3$, see \cite{GL}, \cite{LV}, \cite{Me} and \cite{Va$_2$}.

\medskip

Finally, applying Theorem \ref{th7.8.2.1} in $\Omega$, we come to
the inequality (\ref{eq_1.14.111.1}). $\Box$
\end{proof}

\medskip

The following theorem is also a consequence of Theorem
\ref{th7.8.2.1} and its proof is perfectly similar to the proof of
the last corollary and based on the estimates (\ref{eqOSadm.1}) and
$e^{\tau}\ge \frac{\tau^N}{N!}$ for all $N=1,2,\ldots$.

\medskip

\begin{theorem}
{}\label{th7.8.2.1.i}{\sl Let $D$ be a domain in ${\Bbb
R}^{n},$\thinspace \thinspace $n\ge 2,$ and let $f_{j}$,
$j=1,2,\ldots ,$ be a sequence of ho\-me\-o\-mor\-phisms of $D$ into
${\Bbb R}^{n}$ of the class $W^{1,1}_{\mathrm loc}(D)$ converging
locally uniformly to a homeomorphism $f$ of $D$ into ${\Bbb R}^{n}$.
Then, for every continuous nondecreasing convex function $\psi :
\bar{\mathrm I}\to \overline{\mathbb{R}^+}$ and every open set
$\Omega\subseteq D$,
\begin{equation}\label{eq_1.14.i}
    \int\limits_{\Omega} e^{\psi (P_f (x))}\ dm(x)\ \leq\
    \liminf\limits_{j \to \infty} \int\limits_{\Omega} e^{\psi (P_{f_j} (x))}\
    dm(x)\ .
\end{equation}
If in addition $\psi$ is  nonconstant and the right hand side in
(\ref{eq_1.14.i}) is finite, then $f\in W^{1,p}_{\mathrm
loc}(\Omega)$ for all $p\in [1, n)$ and $f$ is differentiable a.e.
in $\Omega$. Moreover, if
\begin{equation}\label{eq_1.14.z.i}
    \limsup\limits_{j \to \infty} \int\limits_{\Omega} e^{\psi (P_{f_j} (x))}\
    dm(x)\ <\ \infty\ ,
\end{equation}} then $\partial_if_{j}\to\partial_if$ weakly in
$L^p_{\mathrm loc}(\Omega)$ as $j\to\infty$ for all $i=1,\ldots n$
and for all $p\in [1, n)$.\end{theorem}

\bigskip

\begin{corollary}
{}\label{cor7.8.2.1.2A}{\sl In particular, for every $\alpha > 0$,
\begin{equation}\label{eq_1.14.iA}
    \int\limits_{\Omega} e^{\alpha\, P_f (x)}\ dm(x)\ \leq\
    \liminf\limits_{j \to \infty} \int\limits_{\Omega} e^{\alpha\, P_{f_j} (x)}\
    dm(x)
\end{equation}
and \begin{equation}\label{eq_1.14.iB}
    \int\limits_{\Omega} e^{\alpha\, K_f (x)}\ dm(x)\ \leq\
    \liminf\limits_{j \to \infty} \int\limits_{\Omega} e^{\alpha\, K_{f_j} (x)}\
    dm(x)\ .
\end{equation} If the right hand side in either
(\ref{eq_1.14.iA}) or (\ref{eq_1.14.iB}) is finite, then $f\in
W^{1,p}_{\mathrm loc}(\Omega)$ for all $p\in [1, n)$ and $f$ is
differentiable a.e. in $\Omega$. Moreover, if
\begin{equation}\label{eq_1.14.iC}
    \limsup\limits_{j \to \infty} \int\limits_{\Omega} e^{\alpha\, P_{f_j} (x)}\
    dm(x)\ <\ \infty
\end{equation}
or \begin{equation}\label{eq_1.14.iD}
    \limsup\limits_{j \to \infty} \int\limits_{\Omega} e^{\alpha\, K_{f_j} (x)}\
    dm(x)\ <\ \infty\ ,
\end{equation} then $\partial_if_{j}\to\partial_if$ weakly in
$L^p_{\mathrm loc}(\Omega)$ as $j\to\infty$ for all $i=1,\ldots n$
and for all $p\in [1, n)$. }\end{corollary}

\bigskip

\begin{remark}\label{remark4.1.P}
The same is true for any function $\psi(t)$ which is a finite
product of the function $\alpha t^{\beta}$ with $\alpha >0$ and
$\beta > 1$, and some of the functions $[\log(A_1+t)]^{\alpha_1}$,
$[\log\log(A_2+t)]^{\alpha_2},\ \ldots \ $, with $\alpha_m$ and
$A_m\in{\mathbb R }$, $m\in{\mathbb N}$, $t\in [T,\infty ]$,
$\psi(t)\equiv\psi(T)$,  $t\in [1,T]$, with a large enough
$T\in{\mathrm I}$. In particular, choosing here $\beta=n-1$ we
obtain various limit inequalities for $K_f$ if $n\ge 3$.
\end{remark}

\bigskip

Now, let us prove the following useful lemma whose prototype for the
plane case can be found in the work \cite{Lo}. We derive it on the
basis of Theorem \ref{th7.8.2.1}.

\medskip
\begin{lemma}\label{lem1.1.1} {\sl Let $D$ be a domain in ${\Bbb
R}^{n},$\thinspace \thinspace $n\ge 2,$ and let $f_{j}$,
$j=1,2,\ldots ,$ be a sequence of ho\-me\-o\-mor\-phisms of $D$ into
${\Bbb R}^{n}$ of the class $W^{1,1}_{\mathrm loc}(D)$ converging
locally uniformly to a mapping $f:D\rightarrow {\Bbb R}^{n}$ which
is differentiable
 a.e. in $D.$ Then on every open set $\Omega\subseteq D$
\begin{equation}\label{eq_1.14.11}
    \int\limits_{\Omega} \Phi (P_f (x))\ \Psi (x)\ dm(x)\ \leq\
    \liminf\limits_{j \to \infty} \int\limits_{\Omega} \Phi (P_{f_j} (x))\ \Psi (x)\  dm(x)
\end{equation}
for every strictly convex function $\Phi :\bar{\mathrm I} \to
\overline{\mathbb{R}^+}$ which is continuous from the left at the
point $T = \sup\limits_{\Phi(t)<\infty}\ t$ and every uniformly
continuous function $\Psi: {\Bbb R}^{n}\to (0,\infty)$ such that
$1/\Psi$ is locally bounded in ${\Bbb R}^{n}$. }\end{lemma}

\medskip

\begin{proof} As in the proof of Theorem \ref{th7.8.2.1}, by Lemma \ref{le7.8.2.1} and countable
additivity of the integral, we may assume that $\Omega$ is bounded.
On such a set by the conditions of the lemma $0<c\le\Psi\le
C<\infty$. With no loss of generality we may assume also that the
right hand side in (\ref{eq_1.14.11}) is finite and, consequently,
the left hand side in (\ref{eq_1.14}) is finite, too.

Let $C(x,h)$ be a cube centered at a point $x\in D$ whose adges of
the length $h$ are oriented along coordinate axes. It follows by the
uniform continuity of the function $\Psi$ that for every
$\varepsilon>0$ there is $\delta(\varepsilon)>0$ such that
$|\Psi(x)-\Psi(x^{\prime})|$ for all $x\in\Omega$ and $x^{\prime}\in
C(x,h)$.

The system of cubes $C (x, h)$,  $x\in \Omega $,  $h< \min (\rho (x,
\partial \Omega) / \sqrt{n}, \delta (\varepsilon))$, forms
a covering of the set $\Omega$ in the sense of Vitali and we can
choose a countable collection of mutually disjoint cubes $C_m= C
(x_m, h_m)\subseteq\Omega$ in the given covering such that $|\Omega
\setminus \cup C_m| =0$, see e.g. Theorem IV.3.1 in \cite{Sa}.

Correspondingly to Theorem \ref{th7.8.2.1}, for $\varepsilon<c$, we
obtain that

$$
\int\limits_{C_m} \Phi (P_f (x))\ \Psi (x)\ dm(x)\ \leq\
$$
$$
\leq\ (\Psi (x_m)-\varepsilon)\ \int\limits_{C_m} \Phi (P_f (x))\
dm(x)\ +\ 2\varepsilon\ \int\limits_{C_m} \Phi (P_f (x))\  dm(x)\
\leq\
$$
$$
\leq\ (\Psi (x_m)-\varepsilon)\ \liminf\limits_{j \to \infty}
\int\limits_{C_m} \Phi (P_{f_j} (x))\  dm(x)\ +\ 2\varepsilon\
\int\limits_{C_m} \Phi (P_f (x))\  dm(x)\ \leq\
$$
$$
\leq\ \liminf\limits_{j \to \infty} \int\limits_{C_m} \Phi (P_{f_j}
(x))\ \Psi (x)\  dm(x)\
    +\ 2\varepsilon\ \int\limits_{C_m} \Phi (P_f (x))\  dm(x)\
$$

From the last inequality by countable additivity of of the integral
and the Fatou lemma we have that
$$
\int\limits_{\Omega} \Phi (P_f (x))\ (\Psi (x)-2\varepsilon)\ dm(x)\
\leq\ \liminf\limits_{j \to \infty} \int\limits_{\Omega} \Phi
(P_{f_j} (x))\ \Psi (x)\  dm(x)\
$$
and, finally, by arbitrariness of $\varepsilon>0$ we obtain the
inequality (\ref{eq_1.14.11}). $\Box$
\end{proof}

\medskip

\begin{remark}\label{remark5.1}
If the limit mapping $f$ is a homeomorphism and $\Phi$ is of special
forms $\Phi(t)=\psi(t^{\alpha})$, $\alpha >(n-1)^2$ or
$\Phi(t)=e^{\psi(t)}$ with a nonconstant continuous nondecreasing
convex function $\psi$, then we may not assume a priori in Lemma
\ref{lem1.1.1} that $f$ is differentiable a.e. and make in addition
the corresponding conclusions on the weak convergences of the first
partial derivatives as in Corollary \ref{cor7.8.2.1.2} and Theorem
\ref{th7.8.2.1.i}. The proofs of these fact are perfectly similar
and hence we omit details.
\end{remark}

\bigskip

Choosing in Lemma \ref{lem1.1.1} $\Psi(x)=1/(1+|x|^2)^n$ we come to
the following theorems on semicontinuity of dilatations in the mean
with respect to the spherical volume in $\overline{{\Bbb R}^n}$.

\bigskip

\begin{theorem}
{}\label{th7.8.2.1.1}{\sl Let $D$ be a domain in ${\Bbb
R}^{n},$\thinspace \thinspace $n\ge 2,$ and let $f_{j}$,
$j=1,2,\ldots ,$ be a sequence of ho\-me\-o\-mor\-phisms of $D$ into
${\Bbb R}^{n}$ of the class $W^{1,1}_{\mathrm loc}(D)$ converging
locally uniformly to a mapping $f:D\rightarrow {\Bbb R}^{n}$ which
is differentiable a.e. in $D.$ Then on every open set
$\Omega\subseteq D$
\begin{equation}\label{eq_1.14.1.1}
    \int\limits_{\Omega} \Phi (P_f (x))\ \frac{dm(x)}{(1+|x|^2)^n}\ \leq\
    \liminf\limits_{j \to \infty} \int\limits_{\Omega} \Phi (P_{f_j} (x))\ \frac{dm(x)}{(1+|x|^2)^n}
\end{equation}
for every strictly convex function $\Phi :\bar{\mathrm I}
 \to \overline{\mathbb{R}^+}$ which is continuous from the left at the
point $T = \sup\limits_{\Phi(t)<\infty}\ t$.}\end{theorem}

\bigskip

\begin{corollary}
{}\label{cor7.8.2.1.1}{\sl Under the hypotheses of Theorem
\ref{th7.8.2.1.1},
\begin{equation}\label{eq_1.14.111.2}
\int\limits_{\Omega} \Phi (K_f (x))\ \frac{dm(x)}{(1+|x|^2)^n}\
\leq\ \liminf\limits_{j \to \infty} \int\limits_{\Omega} \Phi
(K_{f_j} (x))\ \frac{dm(x)}{(1+|x|^2)^n}
\end{equation}
for every continuous  nondecreasing convex function $\Phi
:\bar{\mathrm I} \to \overline{\mathbb{R}^+}$ and every open set
$\Omega\subseteq D$.}
\end{corollary}

\medskip

\begin{corollary}
{}\label{cor7.8.2.1.3}{\sl Let $D$ be a domain in ${\Bbb
R}^{n},$\thinspace \thinspace $n\ge 2,$ and let $f_{j}$,
$j=1,2,\ldots ,$ be a sequence of ho\-me\-o\-mor\-phisms of $D$ into
${\Bbb R}^{n}$  of the class $W^{1,1}_{\mathrm loc}(D)$ converging
locally uniformly to a homeomorphism $f$ of $D$ into ${\Bbb R}^{n}$.
Then, for every $\alpha >n-1$,
\begin{equation}\label{eq_1.14.111.3}
\int\limits_{\Omega} \Phi (K^{\alpha}_f (x))\
\frac{dm(x)}{(1+|x|^2)^n}\ \leq\
    \liminf\limits_{j \to \infty} \int\limits_{\Omega} \Phi (K^{\alpha}_{f_j} (x))\ \frac{dm(x)}{(1+|x|^2)^n}
\end{equation}
for every continuous nondecreasing convex function $\Phi
:\bar{\mathrm I} \to \overline{\mathbb{R}^+}$ and every open set
$\Omega\subseteq D$. Moreover, if $\Phi$ is not constant and the
right hand side of (\ref{eq_1.14.111.3}) is finite, then $f\in
W^{1,p}_{\mathrm loc}(\Omega)$ with $p=n\alpha /(1+\alpha)>n-1$ and
$f$ is differentiable a.e. in $\Omega$. Finally, if in addition
\begin{equation}\label{eq_1.14.111.3.z}
\limsup\limits_{j \to \infty} \int\limits_{\Omega} \Phi
(K^{\alpha}_{f_j} (x))\ \frac{dm(x)}{(1+|x|^2)^n}\ <\ \infty\ ,
\end{equation}
then $\partial_if_{j}\to\partial_if$ weakly in $L^p_{\mathrm
loc}(\Omega)$ as $j\to\infty$ for all $i=1,\ldots n$.}
\end{corollary}

\bigskip

\begin{theorem}
{}\label{th7.8.2.1.k}{\sl Let $D$ be a domain in ${\Bbb
R}^{n},$\thinspace \thinspace $n\ge 2,$ and let $f_{j}$,
$j=1,2,\ldots ,$ be a sequence of ho\-me\-o\-mor\-phisms of $D$ into
${\Bbb R}^{n}$ of the class $W^{1,1}_{\mathrm loc}(D)$ converging
locally uniformly to a homeomorphism $f$ of $D$ into ${\Bbb R}^{n}$.
Then, for every continuous nondecreasing convex function $\psi
:\bar{\mathrm I} \to \overline{\mathbb{R}^+}$ and every open set
$\Omega\subseteq D$,
\begin{equation}\label{eq_1.14.k}
    \int\limits_{\Omega} e^{\psi (P_f (x))}\ \frac{dm(x)}{(1+|x|^2)^n}\ \leq\
    \liminf\limits_{j \to \infty} \int\limits_{\Omega} e^{\psi (P_{f_j} (x))}\
    \ \frac{dm(x)}{(1+|x|^2)^n}\ .
\end{equation}
Moreover, if $\psi$ is nonconstant and the right hand side in
(\ref{eq_1.14.k}) is finite, then $f\in W^{1,p}_{\mathrm
loc}(\Omega)$ for all $p\in [1, n)$ and $f$ is differentiable a.e.
in $\Omega$. Finally, if
\begin{equation}\label{eq_1.14.z.k}
    \limsup\limits_{j \to \infty} \int\limits_{\Omega} e^{\psi (P_{f_j} (x))}\
    \ \frac{dm(x)}{(1+|x|^2)^n}\ <\ \infty\ ,
\end{equation}} then $\partial_if_{j}\to\partial_if$ weakly in
$L^p_{\mathrm loc}(\Omega)$ as $j\to\infty$ for all $i=1,\ldots n$
and for all $p\in [1, n)$.\end{theorem}

\medskip

Remark \ref{remark4.1.P} on examples of $\psi$ is valid in the case
of the spherical volume.

\bigskip

\begin{corollary}
{}\label{cor7.8.2.1.2B}{\sl In particular, for every $\alpha > 0$,
\begin{equation}\label{eq_1.14.iAV}
    \int\limits_{\Omega} e^{\alpha\, P_f (x)}\ \frac{dm(x)}{(1+|x|^2)^n}\  \leq\
    \liminf\limits_{j \to \infty} \int\limits_{\Omega} e^{\alpha\, P_{f_j} (x)}\
    \frac{dm(x)}{(1+|x|^2)^n}
\end{equation}
and \begin{equation}\label{eq_1.14.iBV}
    \int\limits_{\Omega} e^{\alpha\, K_f (x)}\ \frac{dm(x)}{(1+|x|^2)^n}\ \leq\
    \liminf\limits_{j \to \infty} \int\limits_{\Omega} e^{\alpha\, K_{f_j} (x)}\
    \frac{dm(x)}{(1+|x|^2)^n}\ .
\end{equation} If the right hand side in either
(\ref{eq_1.14.iAV}) or (\ref{eq_1.14.iBV}) is finite, then $f\in
W^{1,p}_{\mathrm loc}(\Omega)$ for all $p\in [1, n)$ and $f$ is
differentiable a.e. in $\Omega$. Moreover, if either
\begin{equation}\label{eq_1.14.iCV}
    \limsup\limits_{j \to \infty} \int\limits_{\Omega} e^{\alpha\, P_{f_j} (x)}\
    \frac{dm(x)}{(1+|x|^2)^n}\ <\ \infty
\end{equation}
or \begin{equation}\label{eq_1.14.iDV}
    \limsup\limits_{j \to \infty} \int\limits_{\Omega} e^{\alpha\, K_{f_j} (x)}\
    \frac{dm(x)}{(1+|x|^2)^n}\ <\ \infty\ ,
\end{equation} then $\partial_if_{j}\to\partial_if$ weakly in
$L^p_{\mathrm loc}(\Omega)$ as $j\to\infty$ for all $i=1,\ldots n$
and for all $p\in [1, n)$. }\end{corollary}

\medskip

 \cc
\section{ On Compactness of Sobolev Homeomorphisms }  \label{17}

Recall that a class of mappings is called compact if it is normal
and closed. Combining the above results on normality and closeness,
we obtain the following results on compactness for the classes of
the Sobolev homeomorphisms.

\medskip

Given a domain $D\subset {\Bbb R}^n,$ $n\ge 3,$ a measurable
function $Q:D\rightarrow\bar{\mathrm I}$ and $z_1,z_2\in D,$
$z^{\,\prime}_1, z_2^{\,\prime}\in {\Bbb R}^n,$ $z_1\ne z_2,$
$z_1^{\,\prime}\ne z_2^{\,\prime}$, denote by ${\mathcal {S}}_{Q}$
the family of all homeomorphisms $f$ of $D$ into ${\Bbb R}^n$ in the
Sobolev class $W^{1,1}_{\mathrm loc}$ such that $K^{n-1}_f(x)\le
Q(x)$ a.e. in $D$ and $f(z_1)=z_1^{\,\prime},$
$f(z_2)=z_2^{\,\prime}$. Similarly, given a function
$\Phi:\bar{\mathrm I} \to \overline{\mathbb{R}^+}$ and
$\alpha\geqslant1$, denote by ${\mathcal S}_{\Phi , \alpha}$ the
family of all homeomorphisms $f$ of $D$ into ${\Bbb R}^n$ in the
Sobolev class $W^{1,1}_{\mathrm loc}$ with the same normalization
such that
\begin{equation}\label{eqOS10.36a} \int\limits_D\ \Phi (K^{\alpha}_f(x))\
\frac{dm(x)}{\left(1+|x|^2\right)^n}\ \leqslant\ 1\ .\end{equation}
Finally, given a function $\psi:\bar{\mathrm I} \to
\overline{\mathbb{R}^+}$ and $M\in\mathbb{R}^+$, denote by
${\mathcal S}^{\psi}_M$ the family of all ho\-meo\-mor\-phisms $f$
of $D$ into ${\Bbb R}^n$ in the Sobolev class $W^{1,1}_{\mathrm
loc}$ with the given normalization such that
\begin{equation}\label{eqOS10.36a} \int\limits_D\ e^{\psi (K_f(x))}\
\frac{dm(x)}{\left(1+|x|^2\right)^n}\ \leqslant\ M\ .\end{equation}

\medskip
\begin{lemma}\label{th1}
{\sl Let $Q\in L_{\mathrm loc}^{\gamma}$ for some $\gamma >1$ and
satisfy the condition (\ref{eq6.3.24}). Then the class
$\mathcal{S}_Q$ is compact. }
\end{lemma}

\medskip

\begin{proof} First of all the family $\mathcal{S}_Q$ is normal by Lemma \ref{lem7.3}.
Now, let us prove its closeness. For this goal, consider an
arbitrary sequense $f_j\in\mathcal{S}_Q$, $j=1,2,\ldots$ such that
$f_j\to f$ locally uniformly as $j\to\infty$. Then by Lemma
\ref{lem7.2} $f$ is a homeomorphism with the given normalization. By
Corollary \ref{cor7.8.2.1.2} with $\alpha =\gamma (n-1)> (n-1)$ and
$\Phi(t)\equiv t$ we obtain that $f\in W^{1,p}_{\mathrm loc}$ with
$p=n\alpha/(1+\alpha)>n-1$ and $K^{n-1}_f(x)\le Q(x)$ a.e. in $D$,
i.e., $f\in{\mathcal {S}}_{Q}$. $\Box$
\end{proof}

\bigskip

\begin{theorem}{}\label{th6.6.1.1}{\sl\,
If $Q\in$ FMO$\ \cap\ L_{\mathrm loc}^{\gamma}$, $\gamma >1$, then
the class $\mathcal{S}_{Q}$ is compact.}
\end{theorem}

\medskip

\begin{proof}
Let $x_0\in D.$ We may consider further that $x_0=0\in D.$ Choosing
a positive $\varepsilon_0<\min\left\{{\rm dist\,}\left(0, \,\partial
D\right),\,\, e^{\,-1}\right\},$ we obtain by Lemma 2.1 in
\cite{RS$_3$} that
$$\int\limits_{\varepsilon<|x|<\varepsilon_0}Q(x)\cdot\psi^n(|x|)
 \ dm(x)\,=\,O
\left(\log\log \frac{1}{\varepsilon}\right)$$
%
where $\psi(t)\,=\,\frac {1}{t\,\log{\frac1t}}$. Note that
$I(\varepsilon,
\varepsilon_0)\,:=\,\int\limits_{\varepsilon}^{\varepsilon_0}\psi(t)\,dt\,=\,
\log{\frac{\log{\frac{1}
{\varepsilon}}}{\log{\frac{1}{\varepsilon_0}}}}.$ Thus, the desired
conclusion follows from Lemma \ref{th1}.
\end{proof}$\Box$

\medskip

The following consequences of Theorem \ref{th6.6.1.1} are obtained
correspondingly by Corollary 2.1 and Proposition 2.1 in the paper
\cite{RS$_3$}.

\medskip

\begin{corollary}\label{cor6.6.2.1}{\sl\,
The class \, $\mathcal{S}_{Q}$ is compact if $Q\in L_{\mathrm
loc}^{\gamma}$, $\gamma >1$ and
\begin{equation}\label{eq6.6.3.1}
\overline{\lim\limits_{\varepsilon\rightarrow 0}}\ \
 \dashint_{\mathrm B( x_0 ,\varepsilon)} Q(x)\ \ dm(x) <
 \infty\qquad\qquad\forall\,\,\,x_0
\in D 
\end{equation}}
\end{corollary}

\begin{corollary}\label{cor6.6.4.1} {\sl\, The class\,
$\mathcal{S}_{Q}$ is compact if $Q\in L_{\mathrm loc}^{\gamma}$ for
some $\gamma >1$ and every $x_0 \in D$ is a Lebesgue point of $Q$.}
\end{corollary}

\medskip

\begin{theorem}{}\label{th6.6.5.1}{\sl\, Let $Q\in L_{\mathrm
loc}^{\gamma}$, $\gamma >1$, and satisfy the condition
\begin{equation}\label{eq6.6.3.2}
\int\limits_{0}^{\varepsilon(x_0)}\frac{dr}{rq_{x_0}^{\frac{1}{n-1}}(r)}=
\infty\qquad\forall\, x_0\in D\end{equation} for some
$\varepsilon(x_0)< {\rm dist}\, (x_0,
\partial D)$ where $q_{x_0}(r)$ denotes the average of
$Q(x)$ over the sphere $|x-x_0|=r.$ Then the class $\mathcal{S}_{Q}$
is compact.}
\end{theorem}

\medskip

\begin{proof}
Fix $x_0\in D$  and set
$I=I(\varepsilon, \varepsilon_0)=\int\limits
_{\varepsilon}^{\varepsilon_0}\psi(t)\,dt,$
$\varepsilon\in(0,\varepsilon_0)$, where
$$\psi(t)\quad=\quad \left \{\begin{array}{rr}
1/[tq^{\frac{1}{n-1}}_{x_0}(t)]\ , & \ t\in (\varepsilon,
\varepsilon_0)\ ,
\\ 0\ ,  &  \ t\notin (\varepsilon,
\varepsilon_0)\ .
\end{array} \right.$$
Note that by the Jensen inequality $I(\varepsilon,
\varepsilon_0)<\infty$ for every $\varepsilon\in (0, \varepsilon_0)$
because $Q\in L_{\mathrm loc}^{1}$. On the other hand, in view of
(\ref{eq6.6.3.2}), we obtain that $I(\varepsilon, \varepsilon_*)>0$
for all $\varepsilon\in (0, \varepsilon_*)$ with some
$\varepsilon_*\in (0, \varepsilon_0).$ Simple calculations also show
that
$$\int\limits_{\varepsilon<|x-x_0|<\varepsilon_*} Q(x)\cdot\psi^n(|x-x_0|)\
dm(x)\ =\ \omega_{n-1}\cdot I(\varepsilon, \varepsilon_*)$$ and
$I(\varepsilon,\varepsilon_*)=o\left(I^n(\varepsilon,\varepsilon_*)\right)$
by (\ref{eq6.6.3.2}). Thus, the conclusion of Theorem
\ref{th6.6.5.1} follows by Lemma \ref{th1}.
\end{proof}$\Box$

\medskip
\begin{corollary}\label{cor6.6.7.1}{\sl\, The class $\mathcal{S}_{Q}$
is compact if $Q\in L_{\mathrm loc}^{\gamma}$, $\gamma >1$, has
singularities only of the logarithmic type of the order which is not
more than $n-1$ at every point $x_0\in D$.}
\end{corollary}

\medskip
\begin{theorem}\label{th4.1A.1}{\sl\, The class\,
$\mathcal{S}_{Q}$ is compact if $Q\in L_{\mathrm loc}^{\gamma}$,
$\gamma >1$, and
\begin{equation}\label{eq6.6.3.3}\int\limits_{\varepsilon<|x-x_0|<\varepsilon_0}\frac{Q(x)}{|x-x_0|^n}\,dm(x)=
o\left(\log^n\frac{1}{\varepsilon}\right)\qquad \forall\,x_0\in
D\end{equation}
as $\varepsilon\rightarrow 0$ for some
$\varepsilon_0=\varepsilon(x_0)<{\rm dist\,}(x_0, \partial D).$}
\end{theorem}

\medskip
\begin{proof} The conclusion of Theorem \ref{th4.1A.1} follows from Lemma \ref{th1} by the
choice $\psi(t)=\frac{1}{t}$ in (\ref{eq6.3.24}). \end{proof} $\Box$

\medskip

\begin{theorem}{}\label{th4.2A.1}{\sl\,The class\,
$\mathcal{S}_{Q}$ is compact if, for some $\gamma > 1$,
\begin{equation}\label{eq7.7A.1}
\int\limits_D\Phi\left(Q^{\gamma}(x)\right)\frac{dm(x)}{\left(1+|x|^2\right)^n}\
\le\ M\ <\ \infty
\end{equation}
for a nondecreasing convex function $\Phi:\bar{\mathrm
I}\rightarrow\overline{{\Bbb R}^+}$ such that
\begin{equation}\label{eq7.77.1}
\int\limits_{\delta}^{\infty}
\frac{d\tau}{\tau\left[\Phi^{-1}(\tau)\right]^{\frac{1}{\gamma(n-1)}}}\
=\ \infty
\end{equation}
for some $\delta>\Phi(1).$}
\end{theorem}

\medskip
\begin{proof}
It follows from Theorem \ref{th6.6.5.1} because the conditions
(\ref{eq7.7A.1})--(\ref{eq7.77.1}) imply the condition
(\ref{eq6.6.3.2}) by Theorem 3.1 in \cite{RS$_2$}.
\end{proof} $\Box$

\medskip

\begin{remark}\label{remark6.1.1.q} Note that the condition (\ref{eq7.77.1}) is not only sufficient but
also necessary for the normality and, consequently, for the
compactness of the classes $\mathcal{S}_{Q}$ with $Q$ satisfying the
integral condition (\ref{eq7.7A.1}), see the corresponding example
in \cite{RS$_2$}, Theorem 5.1. Note also, see Proposition 2.3 in
\cite{RS$_2$}, that the condition (\ref{eq7.77.1}) is equivalent to
the following condition
\begin{equation}\label{eqKR4.4.q}
\int\limits_{\delta}^{\infty}\log\,\Phi(t)\
\frac{dt}{t^{\gamma^{\prime}}}\ =\ \infty\end{equation} for all
$\delta>t_0$ where $t_0:=\sup\limits_{\Phi(t)=0}t,$ $t_0=1$ if
$\Phi(1)>0,$ for $\gamma^{\prime}=1+1/\gamma(n-1)$. Note that
$\gamma^{\prime}< n^{\prime}$ where $\frac{1}{n'}+\frac{1}{n}=1$,
i.e., $n^{\prime}=2$ for $n=2$, $n^{\prime}$ is strictly decreasing
in $n$ and $n^{\prime}=n/(n-1)\to1$ and hence $\gamma^{\prime}\to 1$
as $n\to\infty$.
\end{remark}

\bigskip

Finally, let us give criteria for compactness of Sobolev's classes
with no point-wise dominants of dilatations.

\bigskip

\begin{theorem}{}\label{th4.2A.2}{\sl\, The class
$\mathcal{S}_{\Phi, \alpha}$ is compact for all $\alpha >n-1$ and
all continuous nondecreasing convex functions $\Phi:\bar{\mathrm I}
\rightarrow \overline{{\Bbb R}^+}$ such that
\begin{equation}\label{eq7.77.1.i}
\int\limits_{\delta}^{\infty}
\frac{d\tau}{\tau\left[\Phi^{-1}(\tau)\right]^{\frac{1}{\alpha}}}\
=\ \infty
\end{equation}
for some $\delta>\Phi(1).$}
\end{theorem}

\medskip

\begin{proof}
The family $\mathcal{S}_{\Phi, \alpha}$ is normal by Lemma
\ref{lem7.3.z}. Let $f_j$, $j=1,2,\ldots$ be a sequence of
homeomorphisms in $\mathcal{S}_{\Phi, \alpha}$ such that $f_j\to f$
locally uniformly as $j\to\infty$. Then by Lemma \ref{lem7.2} $f$ is
a homeomorphism with the given normalization. Thus, by Corollary
\ref{cor7.8.2.1.2} the class is closed and hence it is compact.
\end{proof} $\Box$

\medskip

\begin{remark}\label{remark6.1.1} Note again that the condition (\ref{eq7.77.1.i}) is not only sufficient but
also necessary for the normality and, consequently, for the
compactness of the classes $\mathcal{S}_{\Phi , {\alpha}}$, see
Theorem 5.1 in \cite{RS$_2$}. Note also that the condition
(\ref{eq7.77.1.i}) is equivalent to the following condition
\begin{equation}\label{eqKR4.4}
\int\limits_{\delta}^{\infty}\log\,\Phi(t)\
\frac{dt}{t^{\alpha^{\prime}}}\ =\ \infty\end{equation} for all
$\delta>t_0$ where $t_0:=\sup\limits_{\Phi(t)=0}t,$ $t_0=1$ if
$\Phi(1)>0,$ for $\alpha^{\prime}=1+1/\alpha$. Note that
$\alpha^{\prime}<n^{\prime}$ where $n^{\prime}=n/(n-1)\to1$ and,
consequently, $\alpha^{\prime}\to 1$ as $n\to\infty$. Finally,
arguing by contradiction it is easy to see that the condition
(\ref{eqKR4.4}) implies that $\Phi$ is strictly convex.
\end{remark}

\bigskip

\begin{theorem}{}\label{th4.2A.2P}{\sl\, The class
$\mathcal{S}^{\psi}_M$ is compact for all continuous nondecreasing
functions $\psi:\bar{\mathrm I} \rightarrow \overline{{\Bbb R}^+}$
such that $e^{\psi(t^{n-1})}$ is convex, ${\psi^m\left(
t^{\frac{1}{n-1}}\right)}$ is convex for some $m\in{\mathbb N}$ on a
segment $[T,\infty ]$ for some $T\in{\mathrm I}$ and
\begin{equation}\label{eq3!!!3P}
\int\limits_{1}^{\infty}\ \psi(\tau)\ \frac{d\tau}{\tau^2}\ =\
\infty\ .
\end{equation} }
\end{theorem}


\begin{remark}\label{remark6.1.1.P} Note that the convexity of the function
$e^{\psi(t^{n-1})}$ is a more weaker condition than the convexity of
$e^{\psi(t)}$ and all the more than the convexity of $\psi(t)$. The
convexity of the function ${\psi\left( t^{\frac{1}{n-1}}\right)}$ is
conversely stronger than the convexity of $\psi(t)$ but the
convexity of ${\psi^m\left( t^{\frac{1}{n-1}}\right)}$ for $m>1$ may
be weaker. Thus, for the case $m=1$ and $T=1$, the first condition
can be omitted because it follows from the second one. However, the
convexity of $e^{\psi(t^{n-1})}$ is a necessary condition for
compactness of the classes $\mathcal{S}^{\psi}_M$ because
$e^{\psi(K_f)}=e^{\psi(P^{n-1}_f)}$, see Remark \ref{remark5.1}.
Moreover, simple calculations show that the condition
(\ref{eq3!!!3P}) for the function $\psi$ is equivalent to the
condition
\begin{equation}\label{eqKR4.4.P}
\int\limits_{1}^{\infty}\log\,\Phi(t)\ \frac{dt}{t^{n^{\prime}}}\ =\
\infty\ , \ \ \ \ \ n^{\prime}=n/(n-1)\ ,\end{equation} for the
function $\Phi(t)=e^{\psi(t^{n-1})}$ because $\Phi(1)\ge 1$. In
tern, (\ref{eqKR4.4.P}) is equivalent to the condition
\begin{equation}\label{eq7.77.1.i.P}
\int\limits_{1}^{\infty}
\frac{d\tau}{\tau\left[\Phi^{-1}(\tau)\right]^{\frac{1}{n-1}}}\ =\
\infty\ ,
\end{equation}
see Proposition 2.3 in \cite{RS$_2$}. Thus, it follows by Theorem
5.1 in \cite{RS$_2$} that the condition (\ref{eq3!!!3P}) is also
necessary.
\end{remark}

\bigskip

\begin{proof}
The family $\mathcal{S}^{\psi}_M$ is normal by Lemma \ref{lem7.3.z}.
Now, consider a sequence of homeomorphisms
$f_j\in\mathcal{S}^{\psi}_M$, $j=1,2,\ldots$, such that $f_j\to f$
locally uniformly.

\medskip

Set $\Psi(t)={\psi^m\left( t^{\frac{1}{n-1}}\right)}$ for
$t\in[T,\infty ]$ and extend the function $\Psi$ to
$\overline{\mathbb{R}^+}$ by the equality $\Psi(t)\equiv \Psi(T)$
for all $t\in[0,T)$. By construction $\Psi$ is continuous,
nondecreasing, not constant, see (\ref{eqKR4.4}), and convex, see
Proposition I.4.8 in \cite{Bou}. Setting
$t_0=\sup\limits_{\Psi(t)=\Psi(0)} t$ and
$T_0=\sup\limits_{\Psi(t)<\infty} t$, we find that $t_0<T_0$ because
$\Psi$ is continuous, nondecreasing and not constant and, taking
$t_*\in(t_0,T_0)$, we have that
$$
\frac{\Psi(t)-\Psi(0)}{t}\ \ge \frac{\Psi(t_*)-\Psi(0)}{t_*}\ >\ 0\
\ \ \ \ \ \ \ \ \ \  \forall\ t\in [\ t_*\ ,\ \infty\ )
$$
by convexity of $\Psi$, see Proposition I.4.5 in \cite{Bou}, i.e.,
$\Psi(t)\ge at$ for $t\ge t_*$ where
$a={[\Psi(t_*)-\Psi(0)]}/{t_*}>0$. After the replacement
$\tau=t^{\frac{1}{n-1}}$, we have that $\psi^m(\tau)\ge a\tau^{n-1}$
for all $\tau\in[\tau_*,\infty ]$ where
$\tau_*=t^{\frac{1}{n-1}}_*$. Thus, applying the inequality
$e^{\psi}\ge\psi^k/k!$ for all $k=lm$, $l\in{\mathbb N}$, we obtain
by (\ref{eqOS10.36a}) that
$$
\int\limits_{D} K^{l(n-1)}_{f_j} (x)\
\frac{dm(x)}{\left(1+|x|^2\right)^n}\ \le\ \tau_*^lS(D)\ +\
(lm)!M/a^l\ <\ \infty\ \ \ \ \ \ \ \  \forall\ j, l\in\ \mathbb N
$$
where $S(D)$ denotes the spherical volume of $D$. Hence by
Proposition \ref{prOS2.4} $f_j\in W^{1,p}_{\mathrm loc}(D)$ for all
$p\in[1,n)$ and, in particular, $f_j\in W^{1,\varphi}_{\mathrm
loc}(D)$, $j=1,2,\ldots $, where $\varphi(t)=t^{p_*}$ for some
$p_*\in(n-1, n)$. Moreover, since $e^{\psi}\ge\psi^m/m!$, we obtain
by the condition (\ref{eqOS10.36a}) that
\begin{equation}\label{eqOS10.36a.v} \int\limits_D\ {\Psi (K^{n-1}_{f_j}(x))}\
\frac{dm(x)}{\left(1+|x|^2\right)^n}\ \leqslant\ m!M\ <\ \infty\ \ \
\ \ \ \ \  \forall\ j\in\ \mathbb N\ .\end{equation} Thus, by
Theorem \ref{th4.2.v} $f$ is a homeomorphism.

Consequently, by Corollary \ref{cor7.8.2.1.3} and the above estimate
for $K_{f_j}$ we have that
$$
\int\limits_{D} K^{l(n-1)}_{f} (x)\
\frac{dm(x)}{\left(1+|x|^2\right)^n}\ \le\ \tau_*^lS(D)\ +\
(lm)!M/a^l\ <\ \infty\ \ \ \ \ \ \ \  \forall\ l\in\ \mathbb N\ ,
$$
$f\in W^{1,p}_{\mathrm loc}(D)$ for all $p\in[1,n)$ and $f$ is
differentiable a.e. in $D$.

Finally, by Theorem \ref{th7.8.2.1.1} applied to
$\Phi(t)=e^{\psi(t^{n-1})}$ we conlude that $f$ satisfies the
condition (\ref{eqOS10.36a}) and hence $f\in\mathcal{S}^{\psi}_M$.

\medskip

Thus, the class $\mathcal{S}^{\psi}_M$ is closed and hence it is
compact.
\end{proof} $\Box$

\medskip

\begin{remark}\label{remark4.1PP} As examples of convex functions $\psi$
satisfying the condition (\ref{eq3!!!3P}) in Theorem \ref{th4.2A.2P}
can be finite products of the function $\alpha t^{\beta}$, $\alpha
>0$, $\beta \ge 1$, and some of the functions
$[\log(A_1+t)]^{\alpha_1}$, $[\log\log(A_2+t)]^{\alpha_2},\ \ldots \
$, $\alpha_m\ge -1$, $A_m\in{\mathbb R }$, $m\in{\mathbb N}$, $t\in
[T,\infty ]$, $\psi(t)\equiv\psi(T)$,  $t\in [1,T]$, with a large
enough $T\in{\mathrm I}$.
\end{remark}

\medskip

 \cc
\section{ On Precision of Conditions}  \label{17}

As it was noted in Remark \ref{remark6.1.1}, the condition
(\ref{eq7.77.1.i}) in Theorem \ref{th4.2A.2} is necessary and, in
tern, it implies the condition (\ref{eq_1.16}). In this section, we
show that other main conditions on the function $\Phi$ in Theorems
\ref{th7.8.2.1}, \ref{th7.8.2.1.1}, \ref{th4.2A.2} and
\ref{th4.2A.2P} are also necessary. Namely, we prove here the
following result.

\bigskip

\begin{theorem}{}\label{th4.2A.2PPP}{\sl\, Let $D$ be a domain in ${\Bbb
R}^{n},$\thinspace \thinspace $n\ge 2,$ and let $\Phi:\bar{\mathrm
I} \rightarrow \overline{{\Bbb R}^+}$ be a nonconstant function with
$\Phi(\infty)=\infty$. If $\Phi$ is either not convex or not
nondecreasing or not continuous from the left at the point $T =
\sup\limits_{\Phi(t)<\infty}\ t$, then there is a sequence of
ho\-me\-o\-mor\-phisms $f_j$ of $D$ into ${\Bbb R}^{n}$ in the class
$W^{1,1}_{\mathrm loc}(D)$ converging locally uniformly to a mapping
$f:D\rightarrow {\Bbb R}^{n}$ differentiable a.e. such that on every
open bounded set $\Omega\subseteq D$
\begin{equation}\label{eq_1.14}
    \int\limits_{\Omega} \Phi (P_f (x))\, \Psi(x)\, dm(x)\ >\
    \limsup\limits_{j \to \infty} \int\limits_{\Omega} \Phi (P_{f_j}
    (x))\,
\Psi(x)\,    dm(x)
\end{equation}
for every continuous function $\Psi : {\Bbb R}^{n}\to (0,\infty)$
and, in particular, for the functions $\Psi(x)\equiv 1$ and
$\Psi(x)=1/(1+|x|^2)^n$.}
\end{theorem}

\bigskip

The proof of Theorem \ref{th4.2A.2PPP} is based on a series of
lemmas that can have a self--independent interest, too. Their
prototypes in the plane case can be found in the papers
\cite{GR$_1$} and \cite{R$_1$} and in the monograph \cite{GR^*}, see
Lemma 12.1.

\medskip

\begin{lemma}\label{precission1}
{\sl Let $t_1$ and $t_2\in\mathrm I$ and $\lambda \in [0,1]$. Then
there is a sequence of $K$-quasiconformal mappings $f_j : {\Bbb
R}^{n}\to{\Bbb R}^{n}$ with $K=\max\ (\ t_1^{n-1},\ t_2^{n-1})$ and
the dilatations $P_{f_j}$, $j=1,2,\ldots $ taking only two values
$t_1$ and $t_2$ a.e. that converges uniformly in ${\Bbb R}^{n}$ to
an affine mapping $f$ with the dilatation
\begin{equation}\label{eq_2.8}
    P_{f}(x)\ \equiv\ t_0\ \colon =\ \lambda\ t_1 + (1- \lambda)\
    t_2  \ ,
\end{equation}
where the choice of the sequence of measurable sets
\begin{equation}\label{eq_2.10}
    E_j = \{ x \in {\Bbb R}^{n} : P_{f_j} (x) = t_1 \}
\end{equation}
depends only on $\lambda$ but not on $t_1$ and $t_2$. Moreover, for
every function $\Phi : \mathrm I \to \mathbb{R}$ and every
continuous function $\Psi : {\Bbb R}^{n}\to {\Bbb R}^+$, there is
\begin{eqnarray}\label{eq_2.9}
  \lim\limits_{j \to \infty} \int\limits_{{E}} \Phi (P_{f_j}(x))\ \Psi(x) \ dm(x)\ &=&  \\
   &=&
  \{ \lambda \Phi (t_1) + (1- \lambda) \Phi (t_2) \} \int\limits_{{E}}  \Psi(x) \
  dm(x) \nonumber
\end{eqnarray}
on every set $E$ in ${\Bbb R}^{n}$ of a finite Lebesgue measure. }
\end{lemma}

\medskip

\begin{proof}
First, let us give the affine mappings $g_l  : {\Bbb R}^{n}\to{\Bbb
R}^{n}$, $l=0,1,2$, in terms of coordinates $x=(x_1,\ldots , x_n)$
and coordinate functions $y=(y_1,\ldots , y_n)$ setting $y_1\equiv
x_1\ldots , y_{n-1}\equiv x_{n-1}$ and $y_n=t_l\, x_n$, $l=0,1,2$.
It is obvious that $P_{g_l}(x)\equiv t_l$, $l=0,1,2$.

Now, let us fix $j=1, \, 2, \, \dots$, and split the space ${\Bbb
R}^{n}$ into fibers by the hyperplanes $\mathrm  H_{jm}=\{ x\in{\Bbb
R}^{n}: x_n=m\, 2^{-j}\}$, $m=0, \, \pm 1, \, \pm 2, \, \dots$. In
tern, split every such fiber into two fibers by the hyperplane
$\Tilde {\mathrm H}_{jm}=\{ x\in{\Bbb R}^{n}:x_n=m\,
2^{-j}+\lambda\, 2^{-j}\}$.

Next, set $f_j (x)=g_1 (x)$ for all $x\in{\Bbb R}^{n}$ lying between
the hyperplanes $\mathrm  H_{j\, 0}$ and $\Tilde {\mathrm H}_{j\,
0}$. For all the rest fibers, set $f_j (x)=g_l (x) + c^{(l)}_{j\, m
}e_n$, where $l=1$ between the hyperplanes $\mathrm H_{j\, m}$ and
$\Tilde {\mathrm H}_{j\, m}$ and $l=2$ between the hyperplanes
$\Tilde {\mathrm H}_{j\, m}$ and $\mathrm  H_{j\, (m+1)}$ and the
constants $c^{(l)}_{j\, m}$ are found by induction on $m=0, \, \pm
1, \, \pm 2, \, \dots$ on both sides of zero from the gluing
condition. Set also $f(x)\equiv g_0(x)$ in ${\Bbb R}^{n}$.

Note that by the construction $\Delta f_j = \{ \lambda t_1 + (1-
\lambda)t_2 \}\, 2^{-j}\, e_n = t_0\, 2^{-j}\, e_n = \Delta g_0$ for
$\Delta x = 2^{-j} e_n $, see also (\ref{eq_2.8}), and $f_j (x) =
g_1 (x) = g_0 (x)$ on the hyperplane $\mathrm H_{j\, 0}$. Since
$\mathrm H_{j\, m}$ coincide with $\mathrm H_{(j+1)\, 2m}$, the
collection of the given hyperplanes only grows with the growth of
$j$. Consequently, $f_j (x) = f (x)$ on all hyperplanes $\mathrm
H_{j\, m}$, $ j = 1, \, 2, \, \dots$, $m = 0, \, \pm 1, \, \pm 2,
\dots$ that are everywhere dense in ${\Bbb R}^{n}$ and, moreover,
$|f_j(x)-f(x)|\le t_02^{-j}$ for all $ j = 1, \, 2, \, \dots$ and
for all $x\in{\Bbb R}^{n}$. Thus, $f_j(x)\to f(x)$ as $j\to\infty$
uniformly with respect to $x\in{\Bbb R}^{n}$. The relation
(\ref{eq_2.9}) follows immediately from distribution of measures
between values $t_1$ and $t_2$ in the dilatations $P_{f_j}$, $j=1,
2, \, \dots$. \end{proof} $\Box$

\medskip

\begin{lemma}\label{precission2}
{\sl Let $\tau_0$ and $\tau_*\in\mathrm I$ and $\tau_*>\tau_0$. Then
there is a sequence of $K$-quasiconformal mappings $f_j : {\Bbb
R}^{n}\to{\Bbb R}^{n}$ with $K=\tau_*^{(n-1)^2}$ and the dilatations
$P_{f_j}$, $j=1,2,\ldots $ taking only one value $\tau_*$ a.e. that
converges uniformly in ${\Bbb R}^{n}$ to an affine mapping $f$ with
the dilatation $P_f(x)=\tau_0$.}
\end{lemma}

\medskip

\begin{proof} It is easy to see that the conclusion of Lemma \ref{precission2}
follows from the construction given under the proof of Lemma
\ref{precission1} with the choice of $t_0=\tau_0$,
$t_1=\tau_*^{1-n}$ and $t_2=\tau_*$ and of the number
$$
\lambda\ =\ \frac{\tau_*\ -\ \tau_0}{\tau_*-\tau_*^{1-n}}\ \in\
(0,1)
$$
found from the relation $\tau_0=\lambda\, \tau^{1-n}_*+(1-\lambda)\,
\tau_*$. The relations (\ref{eq_2.8})--(\ref{eq_2.9}) themselves are
not relating to Lemma \ref{precission2}.
\end{proof} $\Box$

\medskip

\begin{lemma}\label{precission3}
{\sl Let $\tau_0$ and $\tau_j\in\mathrm I$, $j=1,2,\ldots$, be such
that $\tau_j\to\infty$ as $j\to\infty$. Then there is a sequence of
quasiconformal mappings $f_j : {\Bbb R}^{n}\to{\Bbb R}^{n}$ with the
dilatations $P_{f_j}$, $j=1,2,\ldots $ taking only two values
$\tau_0$ and $\tau_j$ a.e. that converges uniformly in ${\Bbb
R}^{n}$ to a Lipschitz homeomorphism $f$ with the dilatation $P_f$
taking only two values $\tau_0$ and $\infty$ a.e. Moreover, $f$ is
periodic in the variable $x_n$ with the period $1$, $y_n=\psi(x_n)$,
and $y_i=x_i$, $i=1,\ldots , n-1$. Furthermore, we may assume that
$|\{ x\in C:\ P_f(x)=\infty \}|=\lambda$ for a prescribed
$\lambda\in(0,1)$ where $C=\{ x\in{\Bbb R}^{n}: 0\le x_i\le 1,\
i=1,\ldots , n\}$ is the unit cube in ${\Bbb R}^{n}$. }
\end{lemma}

\medskip

\begin{proof} We construct the sequence $f_j$ in the form
$y^{(j)}_i=x_i$, $i=1,\ldots , n-1$, and $y^{(j)}_n=\psi_j(x_n)$
with $\psi_j$ such that $\psi_j(x_n+1)=\psi_j(x_n)+c_j$,
$\psi_j(0)=0$ and, correspondingly, with the constant
$c_j=\psi_j(1)$ by the gluing condition.  Thus, it is sufficient to
define $\psi_j$ on the unit segment $[0,1]$. On the segment $[0,1]$
we apply the procedure for construction of a Cantor set. Namely,
take an arbitrary $\lambda\in(0,1)$ and set $q=1-\lambda^{-1}$.

Let us through out the central interval of the length $q$ from the
segment $[0,1]$ and denote the rest set by $E_1$. Now, form $\psi_1$
applying the contraction in $\tau_1^{n-1}$ times and shifts on the
segments of the set $E_1$ and the tension in $\tau_0$ times and a
shift on $[0,1]\setminus E_1$ and the gluing condition. Arguing by
induction, we through out from each segment of the set $E_j$ central
intervals of equal lengths with the total length $q^{j+1}$ and
denote by $E_{j+1}$ the rest of $E_j$. Then we construct the
function $\psi_j$ on $[0,1]$ applying the contraction in
$\tau_j^{n-1}$ times and shifts on the segments of the set $E_{j+1}$
and  the tension in $\tau_0$ times and  shifts on intervals of the
open set $[0,1]\setminus E_{j+1}$ and the gluing condition. It is
easy to verify that the Cantor set $E=\cap E_j$ has the length
$\lambda$.

Let $E_0$ be the periodic extension of $E$ onto $\Bbb R$. Then it is
clear that $\lim\limits_{j\to\infty}  f_j=f$ where $f$ has the form
$y_i=x_i$, $i=1,\ldots , n-1$, and $y_n=\psi(x_n)$ with $\psi$ such
that $\psi(x_n+1)=\psi(x_n)+c$, $\psi(0)=0$ and the constant
$c=(1-\lambda)\tau_0$ is equal to the length of the set
$[0,1]\setminus E$ multiplied on $\tau_0$, and
$|\psi(x_n^{(1)})-\psi(x_n^{(2)})|\le\tau_0 |x_n^{(1)}-x_n^{(2)}|$
is equal to the total length of all inervals of the open set
$[0,1]\setminus E$ lying in the interval $(x_n^{(1)},x_n^{(2)})$
multiplied on $\tau_0$. Thus, $f$ is a Lipschitz homeomorphism and
hence it is differentiable a.e. Moreover, note that a.e. point of
the set $E_0$ is its density point, see e.g. \cite{Sa}, and,
consequently, we have that $\partial_n f=0$ for a.e. $x_n\in E_0$.
Hence $P_f(x)=\infty$ on a subset of $C$ of the measure $\lambda$.
The rest statements of the lemma are obvious from the construction.
\end{proof} $\Box$

\bigskip

{\it Proof of Theorem \ref{th4.2A.2PPP}.} 1) Let us first assume
that the function $\Phi$ is not convex on $\mathrm I$, i.e., there
exist $t_l\in\mathrm I$, $l=1,2$, $t_1<t_2$ and $\lambda\in(0,1)$
such that
$$
\Phi(\lambda\, t_1+(1-\lambda)\, t_2))\ >\ \lambda\, \Phi(t_1)\ +\
(1-\lambda)\, \Phi(t_2)\ .
$$
However, then by Lemma \ref{precission1} there is a sequence of
homeomorphisms $f_j : {\Bbb R}^{n}\to{\Bbb R}^{n}$ in the class
$W^{1,1}_{\mathrm loc}$ with the dilatations $P_{f_j}$,
$j=1,2,\ldots $ taking only two values $t_1$ and $t_2$ a.e. that
converges uniformly in ${\Bbb R}^{n}$ to an affine mapping $f$ with
the dilatation $P_{f}(x)\ \equiv\ t_0$ where $t_0 = \lambda\, t_1 +
(1- \lambda)\, t_2$ and (\ref{eq_2.9}) holds for every continuous
function $\Psi : {\Bbb R}^{n}\to (0,\infty)$, in particular, for the
functions $\Psi(x)\equiv 1$ and $\Psi(x)=1/(1+|x|^2)^n$, i.e., by
the above assumption
$$
 \lim\limits_{j \to \infty} \int\limits_{{\Omega}} \Phi (P_{f_j}(x))\ \Psi(x) \
 dm(x)\ <\
  \int\limits_{{\Omega}} \Phi (P_{f}(x))\ \Psi(x) \
  dm(x)
$$
for every open bounded set $\Omega$ in $D$.

\medskip

2) Now, let $\Phi$ be not nondecreasing on $\mathrm I$, i.e., there
exist points $\tau_0$ and $\tau_*\in\mathrm I$, $\tau_0<\tau_*$ such
that $\Phi(\tau_0)>\Phi(\tau_*)$. Then by Lemma \ref{precission2}
there is a sequence of homeomorphisms $f_j : {\Bbb R}^{n}\to{\Bbb
R}^{n}$ in the class $W^{1,1}_{\mathrm loc}$ with the dilatations
$P_{f_j}(x)=\tau_*$ a.e., $j=1,2,\ldots $ a.e. that converges
uniformly in ${\Bbb R}^{n}$ to an affine mapping $f$ with the
dilatation $P_f(x)\equiv\tau_0$. Thus, the inequality
(\ref{eq_1.14}) holds for such a sequence.

\medskip

3) Next, let $\Phi$ be not continuous from the left at the point
$\mathrm T = \sup\limits_{\Phi(t)<\infty}\ t<\infty$. Then by the
last point of the proof there is the limit $\Phi(\mathrm T-0)\
\colon =\lim\limits_{t\to\mathrm T-0} \Phi(t)< \Phi(\mathrm T)$. Let
$t_j$, $j=1,2,\ldots $, be a sequence of numbers in $[1,\mathrm T)$
converging to $\mathrm T$, set $t_0=\mathrm T$ and consider the
sequence of the homeomorphisms $f_j  : {\Bbb R}^{n}\to{\Bbb R}^{n}$,
$j=0,1,2,\ldots $, given in terms of coordinates $x=(x_1,\ldots ,
x_n)$ and coordinate functions $y=(y_1,\ldots , y_n):\ $  $y_1\equiv
x_1\ldots , y_{n-1}\equiv x_{n-1}$ and $y_n=t_j\, x_n$,
$j=0,1,2,\ldots $. It is obvious that $f_j\to f\ \colon ={f_0}(x)$
and $P_{f_j}(x)\to P_f(x)$ uniformly in ${\Bbb R}^{n}$ as
$j\to\infty$ because $P_{f_j}(x)\equiv t_j$ and, thus, the
inequality (\ref{eq_1.14}) holds for such a sequence.

Finally, let $\sup\limits_{\Phi(t)<\infty}\ t=\infty$ and let $\Phi$
be not continuous from the left at $\infty$. By the points 1) and 2)
of the proof we may assume that $\Phi$ is nondecreasing and convex
on $\mathrm I$. We may also assume that the function $\Phi$ is
extended from $\bar{\mathrm I}$ to $\overline{\mathbb{R}^+}$ by the
equality $\Phi(t)\equiv \Phi(1)$ for all $t\in[0,1)$. Such extended
function $\Phi$ is nondecreasing and convex, see e.g. Proposition
I.4.8 in \cite{Bou}. If $\Phi$ is not constant on $\mathrm I$, then
$t_0=\sup\limits_{\Phi(t)=\Phi(0)} t < \infty$ and taking
$t_*\in(t_0,\infty )$, we have that
$$
\frac{\Phi(t)-\Phi(0)}{t}\ \ge \frac{\Phi(t_*)-\Phi(0)}{t_*}\ >\ 0\
\ \ \ \ \ \ \ \ \ \  \forall\ t\in [\ t_*\ ,\ \infty\ )
$$
by convexity of $\Phi$, see e.g. Proposition I.4.5 in \cite{Bou},
i.e., $\Phi(t)\ge at$ for $t\ge t_*$ where
$a={[\Phi(t_*)-\Phi(0)]}/{t_*}>0$. However, then $\Phi(t)\to\infty$
as $t\to\infty$, i.e., $\Phi$ is continuous at $\infty$. Thus, it
remains to consider the case when $\Phi(t)$ is constant on $\mathrm
I$. But in this case the conclusion of Theorem \ref{th4.2A.2PPP}
follows from Lemma \ref{precission3}. $\Box$

\bigskip

Similar results on the convergence and compactness  can be obtained
for more general Orlicz--Sobolev classes, however, this requests an
additional research, see e.g. \cite{KRSS}, that will be published
el\-se\-whe\-re.

\medskip
\noindent
{\bf Vladimir Ryazanov, Ruslan Salimov and Evgeny Sevost'yanov:}\\
Institute of Applied Mathematics and Mechanics, NAS of Ukraine,\\
74 Roze Luxemburg Str., Donetsk, 83114, Ukraine,\\
vlryazanov1@rambler.ru, salimov07@rambler.ru, brusin2006@rambler.ru

\end{document}